\documentstyle[11pt,amsmath,amssymb]{article}

\begin{document}

\newtheorem{lem}{Lemma}[section]
\newtheorem{th}{Theorem}[section]
\newtheorem{prop}{Proposition}[section]
\newtheorem{rem}{Remark}[section]
\newtheorem{define}{Definition}[section]
\newtheorem{cor}{Corollary}[section]

\allowdisplaybreaks

\makeatletter\@addtoreset{equation}{section}\makeatother
\def\theequation{\arabic{section}.\arabic{equation}}

\newcommand{\D}{{\cal D}}

\newcommand{\N}{{\Bbb N}}
\newcommand{\C}{{\Bbb C}}
\newcommand{\Z}{{\Bbb Z}}
\newcommand{\R}{{{\Bbb R}^d}}
\newcommand{\Rp}{{\R_+}}
\newcommand{\eps}{\varepsilon}
\newcommand{\om}{\omega}

\newcommand{\supp}{\operatorname{supp}}
\newcommand{\la}{\langle}
\newcommand{\ra}{\rangle}
\newcommand{\const}{\operatorname{const}}

\newcommand{\tob}{\Subset}
\newcommand{\dd}{\overset{{.}{.}}}

\renewcommand{\emptyset}{\varnothing}
\renewcommand{\tilde}{\widetilde}
\newcommand{\rom}[1]{{\rm #1}}
\newcommand{\FC}{{\cal F}C_{\mathrm b}^\infty({\cal D},\Gamma)}
\newcommand{\FCo}{{\cal F}C_{\mathrm b}^\infty({\cal D},\dd\Gamma)}
\newcommand{\FCC}{{\cal F}C_{\mathrm b}^\infty(\overline{{\cal D}},\dd\Gamma)}
\newcommand{\FCD}{{\cal F}C_{\mathrm b}^\infty({\cal D},{\cal D}')}

\newcommand{\vph}{\varphi}
\newcommand{\fii}{\vph}
\newcommand{\di}{\partial}
\renewcommand{\div}{\operatorname{div}}
\newcommand{\Oc}{{\cal O}_{\mathrm c}(\R)}
\newcommand{\Gamman}{\Gamma_0}

\newcommand{\EE}{{\cal E}}

\newcommand{\So}{S_{{\mathrm out},\,\epsilon}}

\begin{center}{\Large \bf Infinite interacting diffusion particles I:\\ Equilibrium process
and its scaling limit\\[2mm]\large  Yuri Kondratiev, Eugene
Lytvynov, and Michael R\"ockner}\end{center}

\begin{abstract}

\noindent A stochastic dynamics $({\bf X}(t))_{t\ge0}$ of a
classical continuous system is a stochastic process which takes
values in the space $\Gamma$ of all locally finite subsets
(configurations) in $\R$ and which has a Gibbs measure $\mu$ as an
invariant measure. We assume that $\mu$ corresponds to a symmetric
pair potential $\phi(x-y)$. An important class of stochastic
dynamics of a classical continuous system is formed by diffusions.
Till now, only one type of such dynamics---the so-called gradient
stochastic dynamics, or interacting Brownian particles---has been
investigated. By using the theory of Dirichlet forms from
\cite{MR92}, we construct and investigate a new type of stochastic
dynamics, which we call infinite interacting diffusion particles.
We introduce a Dirichlet form ${\cal E}_\mu^\Gamma$ on
$L^2(\Gamma;\mu)$, and under general conditions on the potential
$\phi$,  prove its closability. For a potential $\phi$ having a
``weak'' singularity at zero, we also write down an explicit form
of the generator of ${\cal E}_\mu^\Gamma$ on the set of smooth
cylinder functions. We then show that, for any Dirichlet form
${\cal E}_\mu^\Gamma$, there exists a diffusion process that is
properly associated with it. Finally, in a way parallel to
\cite{GKLR}, we study a scaling limit of interacting diffusions in
terms of convergence of the corresponding Dirichlet forms, and we
also show that these scaled processes are tight in
$C([0,\infty),{\cal D}')$, where ${\cal D}'$ is the dual space of
${\cal D}{:=}C_0^\infty(\R)$.

\end{abstract}

\noindent 2000 {\it AMS Mathematics Subject Classification}.
Primary: 60K35, 60B12. Secondary: 60H15, 82C22.\vspace{3mm}

\section{Introduction}

%%%%%%%%%%%%%%%%%%%%%%%%
A stochastic  dynamics $({\bf X}_t)_{t\ge0}$ of a classical
continuous system is a stochastic  process which takes values in
the space $\Gamma$ of all locally finite subsets (configurations)
in $\R$ and which has a Gibbs measure $\mu$ as an invariant
measure. We assume that $\mu$ corresponds to a symmetric,
translation invariant pair potential $\phi(x-y)$ and activity
$z>0$.

An important class of stochastic dynamics of a classical
continuous system is formed by diffusions. Till now, only one type
of such dynamics---the so-called gradient stochastic dynamics, or
interacting Brownian particles---has been investigated.  This
diffusion process  informally solves the following system of
stochastic differential equations:
\begin{align}\label{erdrj} dx(t)&=-\sum_{y(t)\in{\bf X}(t),\, y(t)\ne
x(t)}\nabla\phi(x(t)-y(t))\, dt+\sqrt2\, dB^x(t),\qquad
x(t)\in{\bf X}(t),\\ {\bf X}(0)&=\gamma\in\Gamma,\notag\end{align}
where $(B^x)_{x\in\gamma}$ is a sequence of independent Brownian
motions. The study of such diffusions has been initiated by
R.~Lang \cite{Lang77} (see also \cite{Shi79, Fritz}), who
considered the case $\phi\in C_0^3(\R)$ using finite-dimensional
approximations of stochastic differential equations. More singular
$\phi$, which are of particular interest from the point of view of
statistical mechanics, have been treated by H.~Osada \cite{Osa96}
and M.~Yoshida \cite{Yos96}. These authors were the first to use
the Dirichlet form approach from \cite{MR92} for the construction
of such processes. However, they could not write down the
corresponding generator explicitly, hence could not prove that
their processes actually solve \eqref{erdrj} weakly. This,
however, was proved in \cite{AKR4} (see also the survey paper
\cite{Roeckner}) by showing an integration by parts formula for
the respective Gibbs measures.

But the gradient stochastic dynamics is, of course,  not the
unique diffusion process which has $\mu$ as an invariant measure.
Indeed, let us consider the following system of stochastic
differential equations:
\begin{align}\label{jigvfuzg} dx(t)&=\sqrt2\,
\exp\left[\frac12\sum_{y(t)\in{\bf X}(t),\, y(t)\ne
x(t)}\phi(x(t)-y(t))\right]\,  dB^x(t),\qquad x(t)\in{\bf X}(t),\\
{\bf X}(0)&=\gamma\in\Gamma.\notag\end{align} At least informally,
one  sees that this dynamics does also leave $\mu$ invariant. Note
that, in system \eqref{erdrj}, the information about the
interaction between particles is concentrated in the drift term,
while in system \eqref{jigvfuzg} the interaction is in the
diffusion coefficient and the drift term is absent. This is why we
shall call a process that (weakly) satisfies \eqref{jigvfuzg}
infinite interacting diffusion particles, or just interacting
diffusions. Note that, if there is no interaction ($\phi=0$), both
processes solving \eqref{erdrj} and \eqref{jigvfuzg} coincide.

In this paper, using the Dirichlet form approach (see
\cite{MR92}), under very wide conditions on $\phi$ (more
pricisely, under (A1), (A2) or (A1), (A3), see below), we
construct an equilibrium process that weakly solves
\eqref{jigvfuzg}. The problem of constructing a solution
 for \eqref{jigvfuzg} which starts from a given
configuration, or a given distribution, is still open. Actually,
this problem may be studied by using the ideas developed for the
Hamiltonian and gradient stochastic dynamics,  see \cite{KKO02a},
that is, by obtaining   equations for the time evolution of
correlation functions, or corresponding generating (Bogoliubov)
functionals. This will be the subject of our forthcoming research.

There also exists another type of stochastic dynamics of a
classical continuous system---the so-called  Glauber-type
dynamics, which is a spatial birth-and-death process, see
\cite{KL}.

Let us  briefly describe the contents of the paper. After some
preliminary information about Gibbs measures in
Section~\ref{Section2}, we construct a bilinear form ${\cal
E}_\mu^\Gamma$ on $L^2(\Gamma;\mu)$  in Section~\ref{Section3}.
This form is defined on the set of smooth cylinder functions as
follows:
\begin{equation}\label{jhuig}{\cal E}_\mu^\Gamma(F,G)=\int_\Gamma
\mu(d\gamma) \int_{\R}zm(dx)\, \langle\nabla_x
F(\gamma\cup\{x\}),\nabla_x
G(\gamma\cup\{x\})\rangle,\end{equation} where $m$ denotes
Lebesgue measure on $\R$. By using the Georgii--Nguyen--Zessin
identity, one gets an equivalent representation of ${\cal
E}_\mu^\Gamma$:
\begin{equation}\label{uifgzuzg}
{\cal E}_\mu^\Gamma(F,G)=\int_\Gamma \mu(d\gamma)
\sum_{x\in\gamma}\exp\left[\sum_{y\in\gamma\setminus\{x\}}\phi(x-y)\right]\langle\nabla_x
F(\gamma),\nabla_x G(\gamma)\rangle,\end{equation} where $\nabla_x
F(\gamma)$ is defined as in \eqref{neww}. We show that ${\cal
E}_\mu^\Gamma$ is a pre-Dirichlet form. To compare our situation
with  the gradient dynamics, let us recall that, in the latter
case, the corresponding Dirichlet form $\tilde{\cal E}_\mu^\Gamma$
for $F,G$ as in \eqref{jhuig} looks as follows (see \cite{AKR4}):
\begin{align*}&\tilde{\cal E}_\mu^\Gamma(F,G)= \int_\Gamma \mu(d
\gamma)\sum_{x\in\gamma}\langle\nabla_x F(\gamma),\nabla_x
G(\gamma)\rangle\\ &\qquad=\int_\Gamma \mu(d\gamma)\int_{\R}z
m(dx)\,
\exp\left[-\sum_{y\in\gamma\setminus\{x\}}\phi(x-y)\right]\langle\nabla_x
F(\gamma\cup\{x\}),\nabla_x G(\gamma\cup\{x\})\ra.\end{align*}

In Section~\ref{ersews6543}, we prove the closability of ${\cal
E}_\mu^\Gamma$ under fairly general conditions on the potential
$\phi$, and in Section~\ref{Section5}, assuming additionally that
the function $e^{\phi(x)}$ is integrable in a neighborhood of zero
(which still admits some ``weak'' singularity of $\phi$), we write
down the generator of ${\cal E}_\mu^\Gamma$ on the set of smooth
cylinder functions. This generator looks as follows:
\begin{equation}\label{ghhsb} H_\mu^\Gamma
F(\gamma)=-\sum_{x\in\gamma}\exp\left[\sum_{y\in\gamma\setminus\{x\}}\phi(x-y)\right]\Delta_x
F(\gamma).\end{equation} In Section~\ref{Section6}, we prove the
existence of a conservative diffusion process which is properly
associated with the (closed) Dirichlet form ${\cal E}_\mu^\Gamma$.
This process lives, in general, in the bigger space $\dd\Gamma$ of
all locally finite multiple configurations, but we prove,
analogously to \cite{RS98}, that the process indeed lives in
$\Gamma$ in  case $d\ge2$. (If $d=1$, one cannot, of course,
exclude collisions of particles.) According to \eqref{ghhsb}, the
constructed diffusion process informally solves \eqref{jigvfuzg}.

 Section~\ref{Section7} is devoted to the study of a scaling
limit of the constructed  process $({\bf X}(t))_{t\ge0}$. The
scaling we study is the same as the one considered by many authors
for the gradient stochastic dynamics. The scaled process $({\bf
X}_\epsilon(t))_{t\ge0}$ is defined by
\begin{equation}\label{hvff}{\bf X}_\epsilon(t){:=}
S_{\mathrm{out},\,\epsilon}(S_{\mathrm{in},\,\epsilon}({\bf
X}(\epsilon^{-2}t))),\qquad t\ge0,\ \epsilon>0,\end{equation} and
we are interested in the scaling limit as $\epsilon\to0$. The
first scaling in \eqref{hvff}, $S_{\mathrm{in},\,\epsilon}$,
scales the positions of  particles inside  the configuration space
as follows: $$ \Gamma\ni\gamma\mapsto S_{\mathrm{in},\,\epsilon}
(\gamma){:=}\{\epsilon x\mid x\in\gamma\}\subset\Gamma.$$ The
second scaling, $S_{\mathrm{out},\,\epsilon}$, leads out of the
configuration space and is given by $$ \Gamma\ni\gamma\mapsto
S_{\mathrm{out},\,\epsilon} (\gamma){:=}
\epsilon^{d/2}\gamma-\epsilon^{-d/2}\rho\, dx\in{\cal D'},$$ where
we identify the configuration with the corresponding sum of Dirac
measures, $\rho$ is the first correlation function of $\mu$, and
${\cal D}'$ is the dual space of ${\cal D}{:=}C_0^\infty(\R)$.

T. Brox showed in \cite{Brox} that, in the low activity-high
temperature regime, the Gibbs measure $\mu$ converges under the
scaling $S_{\mathrm{out},\,\epsilon}S_{\mathrm{in},\,\epsilon}$ to
a corresponding white noise measure $\nu_c$ with covariance
operator $c\operatorname{Id}$, where the constant $c>0$ is
explicitly given through the first and second moments of $\mu$.
However, T.~Brox believed that there is no limiting Markov process
for the scaling limit of the gradient stochastic dynamics. Then,
H.~Rost gave some heuristic arguments for the existence of a
limiting generalized Ornstein--Uhlenbeck process \cite{Rost}. In
the celebrated paper \cite{Spohn}, H. Spohn described a proof of
convergence of the scaled processes in the case where the
underlying potential $\phi$ is smooth, compactly supported, and
positive, and $d\le3$. In \cite{GP}, M.~Z.~Guo and G.~Papanicolaou
tried to prove convergence of the corresponding resolvents,
however their considerations were on a more heuristic level.
Finally, in the recent paper \cite{GKLR}, in the case of a general
potential $\phi$, the authors proved convergence of the processes
on the level of convergence of the associated Dirichlet forms.
Furthermore, the tightness of the processes in $C([0,\infty),{\cal
D }')$ was proven, and the convergence of the processes in law was
shown  under the assumption that the Boltzmann--Gibbs principle
holds.

In this paper, we follow the approach of \cite{GKLR}. So, we show
that, on the level of convergence of the associated Dirichlet
forms, the scaled processes $({\bf X}_\epsilon(t))_{t\ge0}$
converge to the generalized Ornstein--Uhlenbeck process $({\bf
N}(t))_{t\ge0}$ in ${\cal D}'$ that  informally  satisfies the
following stochastic differential equation:
\begin{equation}\label{uhisidf} d{\bf N}(t,x)=\frac 1c\,\Delta{\bf
N}(t,x)\, dt+\sqrt2\, d{\bf W}(t,x),
\end{equation} where $({\bf W}(t))_{t\ge0}$ is a Brownian motion on ${\cal
D}'$ with covariance operator $-\Delta$. We recall that the
limiting process $(\tilde{\bf N}(t))_{t\ge0}$ of the gradient
stochastic dynamics satisfies: \begin{equation}\label{jsdgbv}
d\tilde{\bf N}(t,x)=\frac\rho c\, \Delta\tilde{\bf N}(t,x)\,
dt+\sqrt{2\rho}\, d{\bf W}(t,x). \end{equation} Thus, comparing
\eqref{uhisidf} and \eqref{jsdgbv}, we see that the gradient
stochastic dynamics and the interacting diffusions have  great
similarity on the macroscopic level, though they have different
bulk diffusion coefficients: $\rho/c$ for the former stochastic
dynamics, and $1/c$ for the latter.

We finish this paper by proving the tightness of the scaled
processes $({\bf X}_\epsilon(t))_{t\ge0}$ in $C([0,\infty),{\cal D
}')$. To complete the proof of the convergence in law of the
scaled processes, one still needs to prove the Boltzmann--Gibbs
principle in our situation, which remains an open problem.

It is also possible to study an invariance principle (scaling
limit) of a tagged particle of interacting diffusions (cf.\
\cite{osada, osada1}). This will be the subject of future
research.

Finally, we would like to mention that, though  some  proofs of
the results of this paper use the ideas and techniques developed
for the gradient dynamics,  for convenience of the reader we have
tried to make this paper self-contained as possible.

\section{Gibbs measures on configuration spaces}\label{Section2}

The configuration space $\Gamma:=\Gamma_\R$ over $\R$, $d\in\N$,
is defined as the set of all subsets of $\R$ which are locally
finite: $$\Gamma_\R:=\big\{\,\gamma\subset \R\mid
|\gamma_\Lambda|<\infty\text{ for each compact }\Lambda\subset
\R\,\big\},$$ where $|\cdot|$ denotes the cardinality of a set and
$\gamma_\Lambda:= \gamma\cap\Lambda$. One can identify any
$\gamma\in\Gamma$ with the positive Radon measure
$$\sum_{x\in\gamma}\eps_x\in{\cal M}(\R), $$ where  $\eps_x$ is
the Dirac measure with mass at $x$, and  ${\cal M}(\R)$
 stands for the set of all
positive  Radon  measures on the Borel $\sigma$-algebra ${\cal
B}(\R)$. The space $\Gamma$ can be endowed with the relative
topology as a subset of the space ${\cal M}(\R)$ with the vague
topology, i.e., the weakest topology on $\Gamma$ with respect to
which  all maps $$\Gamma\ni\gamma\mapsto\la f,\gamma\ra:=\int_\R
f(x)\,\gamma(dx) =\sum_{x\in\gamma}f(x),\qquad f\in{\cal D},$$ are
continuous. Here, $\D:=C_0^\infty(\R)$ is the space of all
infinitely differentiable real-valued functions on $\R$ with
compact support. We shall denote by ${\cal B}(\Gamma)$ the Borel
$\sigma$-algebra on $\Gamma$.

Let $\pi_z$, $z>0$, denote the Poisson measure on $(\Gamma,{\cal
B}(\Gamma))$ with intensity measure $zm( dx)$.
 This measure can be characterized by its Laplace
transform $$\int_{\Gamma} \exp[\la f,\gamma\ra]\,\pi_z(d\gamma)
=\exp\bigg(\int_\R(e^{f(x)}-1)\,zm(dx)\bigg),\qquad f\in\D.$$ We
refer  e.g.\ to \cite{VGG,AKR3} for a detailed discussion of the
construction of the Poisson measure on the configuration space.

Now, we proceed to consider Gibbs measures. A pair potential is a
Borel  measurable function $\phi\colon \R\to {\Bbb
R}\cup\{+\infty\}$ such that $\phi(-x)=\phi(x)$ for all $x\in\R$.
We shall also suppose that $\phi(x)\in {\Bbb R}$ for all
$x\in\R\setminus\{0\}$. Let $\Oc$ denote the set of all open,
relatively compact sets in $\R$. Then, for $\Lambda\in\Oc$, the
conditional energy $E_\Lambda^\phi\colon\Gamma\to {\Bbb
R}\cup\{+\infty\}$ is defined by \begin{equation}\label{fgftftrtf}
E_\Lambda^\phi(\gamma):=\begin{cases}\sum\limits_{\{x,y\}\subset\gamma,\,
\{x,y\}\cap\Lambda\ne\varnothing }\phi(x-y),&\text{if }
\sum\limits_{\{x,y\}\subset\gamma,\,
\{x,y\}\cap\Lambda\ne\varnothing}|\phi(x-y)|<\infty,\\
+\infty,&\text{otherwise}.\end{cases}\end{equation} Given
$\Lambda\in\Oc$, define for $\gamma\in\Gamma$ and $\Delta\in {\cal
B}(\Gamma)$
\begin{align}
\Pi_\Lambda^{z,\phi}(\gamma,\Delta){:=}&{\bf
1}_{\{Z_\Lambda^{z,\phi}<\infty\}}
(\gamma)\,[Z_\Lambda^{z,\phi}(\gamma)]^{-1}\notag\\ &\times
\int_\Gamma {\bf 1}_\Delta(\gamma_{\Lambda^{\mathrm
c}}+\gamma_\Lambda') \exp\big[
-E_\Lambda^\phi(\gamma_{\Lambda^{\mathrm c }}+\gamma_\Lambda')
\big]\,\pi_z(d\gamma'),\label{awawiuiz}\end{align} where
$\Lambda^{\mathrm c}{:=}\R\setminus\Lambda$ and
\begin{equation}\label{dsesaewa}
Z_\Lambda^{z,\phi}(\gamma){:=}\int_\Gamma\exp\big[
-E_\Lambda^\phi(\gamma_{\Lambda^{\mathrm c}}+\gamma_\Lambda')
\big]\,\pi_z(d\gamma').\end{equation}

A probability measure $\mu $ on $(\Gamma,{\cal B}(\Gamma))$ is
called a grand canonical Gibbs measure with interaction potential
$\phi$ if it satisfies the Dobrushin--Lanford--Ruelle equation
\begin{equation}\label{fdrse}\mu\Pi_\Lambda^{z,\phi}=\mu\qquad\text{for all
}\Lambda\in\Oc.\end{equation} Let ${\cal G}(z,\phi)$ denote the
set of all such probability measures $\mu$.

%It can be shown \cite{Ge79} that the unique grand canonical Gibbs
%measure corresponding to the free case, $\phi=0$, is the Poisson
%measure $\pi_z$.

We rewrite the conditional energy $E_\Lambda^\phi$ in the
following form $$
E_\Lambda^\phi(\gamma)=E_\Lambda^\phi(\gamma_\Lambda)+W(\gamma_\Lambda\mid
\gamma_{\Lambda^{\mathrm c}}),$$ where the term
$$W(\gamma_\Lambda\mid \gamma_{\Lambda^{\mathrm
c}}){:=}\begin{cases}\sum\limits_{x\in\gamma_\Lambda,\,y\in\gamma_
{\Lambda^{\mathrm c}}}\phi(x-y),&\text{if
}\sum\limits_{x\in\gamma_\Lambda,\,y\in\gamma_ {\Lambda^{\mathrm
c}}}|\phi(x-y)|<\infty,\\ +\infty,&\text{otherwise,}\end{cases}$$
describes the interaction energy between $\gamma_\Lambda$ and
$\gamma_{\Lambda^{\mathrm c}}$.  Analogously, we define
$W(\gamma'\mid\gamma'')$ when $\gamma'\cap \gamma''=\varnothing$.

Any $\mu\in{\cal G}(z,\phi)$ satisfies the Georgii--Nguyen--Zessin
identity
\begin{equation}\int_\Gamma \mu(d\gamma)\int_{\R} \gamma(dx) \,
F(\gamma,x) =\int_\Gamma \mu(d\gamma)\int_{\R}
zm(dx)\,\exp\big[-W(\{x\}\mid\gamma )\big] F(\gamma+\varepsilon
_x,x),\label{fdrtsdrt}\end{equation} where
$F:\Gamma\times\R\to[0,+\infty]$ is a measurable function
(\cite[Theorem~2]{NZ}, see also \cite[Theorem~2.2.4]{Kuna}).
 In fact, this
identity uniquely characterizes the Gibbs measures in the sense
that any probability measure $\mu$ on $(\Gamma,{\cal B}(\Gamma))$
belongs to ${\cal G}(z,\phi)$ if and only if $\mu$ satisfies
\eqref{fdrtsdrt}, cf.~\cite[Theorem.~2]{NZ}.

% We suppose that the interaction
%potential $\phi$ is stable, i.e., the following condition is
%satisfied:
%\begin{description}

Let us now describe two classes of Gibbs measures which appear in
classical statistical mechanics of continuous systems
\cite{Ru69,Ru70}. For every $r=(r^1,\dots,r^d)\in\Z^d$, we define
the cube $$Q_r:=\left\{\, x\in\R\mid r^i-\frac 12\le
x^i<r^i+\frac12 \,\right\}.$$ These cubes form a partition of
$\R$. For any $\gamma\in\Gamma$, we set $\gamma_r:=\gamma_{Q_r}$,
$r\in\Z^d$. For $N\in\N$ let $\Lambda_N$ be the cube with side
length $2N-1$ centered at the origin in $\R$, $\Lambda_N$ is then
a union of $(2N-1)^d$ unit cubes of the form $Q_r$.

For $\Lambda\subset\R$, by $\Gamma_\Lambda$ we denote the subset
of $\Gamma$ consisting of all configurations $\gamma\in\Gamma$
such that $\gamma=\gamma_\Lambda$.

 Now, we formulate conditions on
the interaction.

\begin{description}
\item[(S)] ({\it Stability})
 There exists $B\ge0$ such that, for any $\Lambda\in{\cal O}_c(\R)$
and for all $\gamma\in\Gamma_\Lambda$,
$$E_{\Lambda}^\phi(\gamma)\ge -B|\gamma|.$$
\end{description}

Notice that the stability condition automatically implies that the
potential $\phi$ is semi-bounded from below.

\begin{description}

\item[(SS)] ({\it Superstability})
 There exist $A>0$, $B\ge0$ such that, if $\gamma\in\Gamma_
{\Lambda_N}$ for some $N$, then
$$E_{\Lambda_N}^\phi(\gamma)\ge\sum_{r\in\Z^d}\big(A|\gamma_r|
^2-B|\gamma_r|\big).$$

\end{description}

This condition is evidently stronger than (S).

\begin{description}

\item[(LR)] ({\it Lower regularity}) There exists a decreasing positive
function $a\colon\N\to{\Bbb R}_+$ such that
$$\sum_{r\in\Z^d}a(\|r\|)<\infty$$ and for any
$\Lambda',\Lambda''$ which are finite unions of cubes $Q_r$ and
disjoint, with $\gamma'\in\Gamma_{\Lambda'}$,
$\gamma''\in\Gamma_{\Lambda''}$,
$$W(\gamma'\mid\gamma'')\ge-\sum_{r',r''\in\Z^d}a(\|r'-r''\|)
|\gamma_{r'}'|\,|\gamma_{r''}''|.$$ Here, $\|\cdot\|$ denotes the
maximum norm on $\R$.

\item[(I)] ({\it Integrability}) We have
$$\int_\R|1-e^{-\phi(x)}|\,m(dx)<+\infty.$$

\end{description}

We also need

\begin{description}

\item[(UI)] ({\it Uniform integrability}) We have
$$\int_\R|1-e^{-\phi(x)}|\,m(dx)< z^{-1}\exp(-1-2B), $$ where $B$
is as in (S).

\end{description}

A probability measure $\mu$ on $(\Gamma,{\cal B}(\Gamma))$ is
called tempered if $\mu$ is supported by
$$S_\infty{:=}\bigcup_{n=1}^\infty S_n,$$ where $$S_n:=\left\{\,
\gamma\in\Gamma\mid \forall N\in\N\ \sum_{r\in\Lambda_N\cap\Z^d}
|\gamma_r|^2\le n^2|\Lambda_N \cap\Z^d| \,\right\}.$$ By ${\cal
G}^t(z,\phi)\subset{\cal G}(z,\phi)$ we denote the set of all
tempered grand canonical Gibbs measures (Ruelle measures for
short). Due to \cite{Ru70} the set ${\cal G}^t(z,\phi)$ is
non-empty for all $z>0$ and any potential $\phi$ satisfying
conditions (SS), (LR), and (I). Furthermore, the set ${\cal
G}(z,\phi)$ is not empty for potentials satisfying (S) and (UI),
or equivalently, for stable potentials in the low activity-high
temperature  regime, see e.g.\ \cite{MM91, Min67}. A measure
$\mu\in{\cal G}(z,\phi)$ in the latter case is constructed as a
limit of finite volume Gibbs measures corresponding to empty
boundary conditions.

Let us now recall  the so-called Ruelle bound (cf.\ \cite{Ru70}).

\begin{prop}\label{waessedf} Suppose that  either  conditions
\rom{(SS), (LR), (I)} are satisfied and $\mu\in{\cal
G}^t(z,\phi)$\rom, $z>0$\rom, or conditions \rom{(S), (UI)} are
satisfied and $\mu\in{\cal G}(z,\phi)$ is the Gibbs measure
constructed as a limit of finite volume Gibbs measures
 with empty boundary conditions.
Then\rom, for any $n\in\N$\rom, there exists a non-negative
measurable symmetric function $k_\mu^{(n)}$ on $({\Bbb R}^d)^n$
 such
that, for any measurable symmetric function
$f^{(n)}:(\R)^n\to[0,\infty]$,
\begin{align*} &\int_\Gamma \sum_{\{x_1,\dots,x_n\}\subset\gamma}
f^{(n)}(x_1,\dots,x_n)\,\mu(d\gamma)\\&\qquad =\frac1{n!}\,
\int_{(\R)^n} f^{(n)}(x_1,\dots,x_n)
k_\mu^{(n)}(x_1,\dots,x_n)\,m(dx_1)\dotsm m(dx_n),\end{align*} and
\begin{equation}\label{swaswea957}\forall (x_1,\dots,x_n)\in(\R)^n:\quad
k_\mu^{(n)}(x_1,\dots,x_n)\le \xi^n,\end{equation} where $\xi> 0$
is independent of $n$.
\end{prop}

The functions $k_\mu^{(n)}$, $n\in\N$, are called correlation
functions of the measure $\mu$, while \eqref{swaswea957} is called
the Ruelle bound.

The above proposition particularly implies that, for any
$\fii\in{\cal D}$, $\fii\ge0$, and $n\in\N$,
\begin{equation}\label{sweqaw}\int_\Gamma \la \fii,\gamma\ra^n\,\mu(d\gamma)<\infty,\end{equation} that
is, any  measure $\mu$ as in Proposition~\ref{waessedf}  has all
local moments finite.

\section{The bilinear form ${\cal E}_\mu^\Gamma$}\label{Section3}

In what follows, we fix a measure $\mu$ as in
Proposition~\ref{waessedf}. In this section, we shall construct a
 pre-Dirichlet form ${\cal E}_\mu^\Gamma$ on the space $L^2(\Gamma,\mu)$.

 We introduce
the set $\FC$ of all functions   of the form
\begin{equation}\label{3}
\Gamma\ni \gamma\mapsto
F(\gamma)=g_F(\la\fii_1,\gamma\ra,\dots,\la\fii_N,\gamma\ra),
\end{equation}
where $N\in\N$, $\fii_1,\dots,\fii_N\in\D$, and $g_F\in
C^\infty_{\mathrm b}({\Bbb R}^N)$.

 We  define
\begin{equation}\label{drsrs} {\cal E}^\Gamma_\mu(F,G){:=} \int_\Gamma
\mu(d\gamma)\int_{\R}zm( dx)\,\la \nabla_x
F(\gamma+\eps_x),\nabla_x G(\gamma+\eps_x)\ra,\end{equation} where
$F,G\in \FC$. Here, $\nabla_x$  denotes the gradient in the $x$
variable and $\la\cdot,\cdot\ra$ stands for the scalar product in
$\R$. For any $F$ of the form \eqref{3}, we have
\begin{align*}{\cal E}_\mu^\Gamma(F)&\le
\int_\Gamma\mu(d\gamma)\int_{\R}zm(dx)\,\left(\sum_{i=1}^N
|\partial_i
g_F(\la\fii_1,\gamma+\eps_x\ra,\dots,\la\fii_N,\gamma+\eps_x\ra)|
\,|\nabla\varphi_i(x)|\right)^2\\ &\le
\int_\Gamma\mu(d\gamma)\int_{\R}zm(dx)\,\operatorname{const}\sum_{i=1}^N
|\nabla\varphi_i(x)|^2<\infty,\end{align*}  where $\di_j\, g_F$
means derivative with respect to the $j$-th coordinate and,
 as usual, we
set  ${\cal E}^\Gamma_\mu(F){:=}{\cal E}_\mu^\Gamma(F,F)$. Thus,
the the right-hand side of \eqref{drsrs} is well-defined.

In order to get an alternative representation of the form ${\cal
E}_\mu^\Gamma$, we shall suppose the following:

\begin{description}

\item[(A1)] There exists $r>0$ such that $$\sup_{x\in B(r)^{
c}}\phi(x)<\infty,$$ where $B(r)$ denotes the closed ball in $\R$
of radius $r$ centered at the origin.

\end{description}

\begin{lem} \label{sdawawa} In addition to the conditions of
Proposition~\rom{\ref{waessedf},} let  $\phi$ also satisfy
\rom{(A1).} Then\rom, for $\mu\otimes dx$-a\rom.e\rom.\
$(\gamma,x)\in\Gamma\times \R$
\begin{equation}\label{qwqw}\sum_{y\in\gamma}|\phi(x-y)|<\infty\end{equation} and for
$\mu$-a\rom.a\rom.\ $\gamma\in \Gamma$\rom:
\begin{equation}\label{qwqwuqw}\sum_{y\in\gamma\setminus\{x\}}|\phi(x-y)|<\infty\qquad
\text{\rom{for each }}x\in\gamma.\end{equation}\end{lem}

\noindent {\it Proof}. To show \eqref{qwqw}, it suffices to prove
that, for any $\Lambda\in{\cal O}_c(\R)$,
\begin{equation}\label{hgft}
\sum_{y\in\gamma_{(\Lambda^r)^c } } |\phi(x-y)|<\infty\qquad
\text{for $\mu\otimes m$-a.e.\ }(\gamma,x)\in\Gamma\times \Lambda,
\end{equation}
where $\Lambda^r{:=}\{y\in\R: d(y,\Lambda)\le r\}$, $d(y,\Lambda)$
denoting the distance from $y$ to $\Lambda$. By (I) and (A1), we
have \begin{equation}\label{gfuzag} \int_{B(r)^c}|\phi(x)|\,
m(dx)<\infty.\end{equation} Therefore, by
Proposition~\ref{waessedf}
\begin{align*}&\int_\Gamma \mu(d\gamma)\int_\Lambda m(dx)
\sum_{y\in \gamma_{(\Lambda^r)^c}}|\phi(x-y)|\\ &\qquad
=\int_\Lambda m(dx)\int_\Gamma \mu(d\gamma)\int_{\R}\gamma(dy)\,
|\phi(x-y)| {\bf 1} _{(\Lambda^r)^c}(y)\\ &\qquad =\int_\Lambda
m(dx)\int_{\R}m(dy)k_\mu^{(1)}(y)|\phi(x-y)| {\bf 1}
_{(\Lambda^r)^c}(y) \\ &\qquad \le \xi\int_\Lambda m(dx) \int
_{(\Lambda^r)^c} m(dy)\,|\phi(x-y)|\\ &\qquad \le\xi
m(\Lambda)\int_{B(r)^c}|\phi(y)|\, m(dy)<\infty,\end{align*} which
implies \eqref{hgft}.

Analogously, to show \eqref{qwqwuqw} it suffices to prove that,
for any $\Lambda\in{\cal O}_c(\R)$,
\begin{equation}\label{gzhgftf}
\sum_{x\in\gamma_\Lambda}\sum_{y\in
\gamma_{(\Lambda^r)^c}}|\phi(x-y)|<\infty\qquad \text{for
$\mu$-a.e.\ }\gamma\in\Gamma.\end{equation} By
Proposition~\ref{waessedf} and \eqref{gfuzag}
\begin{align*}&\int_\Gamma \sum_{x\in\gamma_\Lambda}\sum_{y\in
\gamma_{(\Lambda^r)^c}}|\phi(x-y)|\, \mu(d\gamma)\\&\qquad
=2\int_{(\R)^2} {\bf 1}_\Lambda (x) {\bf
1}_{(\Lambda^r)^c}(y)|\phi(x-y)|k^{(2)}_\mu(x,y)\,m(dx)\,m(dy)\\&\qquad
\le2\xi^2 m(\Lambda) \int_{B(r)^c}|\phi(y)|\,
m(dy)<\infty,\end{align*} which implies \eqref{gzhgftf} \quad
\vspace{2mm}$\blacksquare$

By using \eqref{fdrtsdrt},  \eqref{drsrs}, and
Lemma~\ref{sdawawa}, we have for any $F,G\in\FC$ \begin{align}
{\cal E}_\mu^\Gamma(F,G)&=\int_{\Gamma}\mu(d\gamma)\int_{\R}z
m(dx)\, \exp\left[ -\sum_{y\in\gamma}\phi(x-y)\right]\notag\\&
\quad\times \exp\left[ \sum_{y\in\gamma}\phi(x-y)\right]
\la\nabla_x F(\gamma+\eps_x),\nabla_xG(\gamma+\eps_x)\ra \notag\\
&= \int_\Gamma S^\Gamma (F,G)\,d\mu.\label{waeaw}
\end{align} Here, \begin{gather} S^\Gamma (F,G) (\gamma){:=}\int_\R
A(\gamma,x) \la \nabla_x F(\gamma),\nabla_x
G(\gamma)\ra\,\gamma(dx),\notag\\
A(\gamma,x){:=}\exp\left[\sum_{y\in\gamma\setminus\{x\}}\phi(x-y)\right],\qquad
 x\in\gamma,\ \text{$\mu$-a.e.\ $\gamma\in\Gamma$}, \label{iouvdc}
\end{gather}
and for any $\gamma\in\Gamma$ and $x\in\R^d$
\begin{equation}\label{neww}\nabla_x F(\gamma){:=}\nabla_y
F(\gamma-\varepsilon_x+\varepsilon_y)\big|_{y=x}.\end{equation}
Note that, since $F\in \FC$, it naturally extends from $\Gamma$ to
all of ${\cal D}'{:=}$dual of $\cal D$. In particular, $F(\gamma)$
is defined if $\gamma$ is a signed measure such that $|\gamma|$ is
finite on compacts.

\begin{lem}\label{kjhdrd} Let the conditions of Lemma~\rom{\ref{sdawawa}} be satisfied\rom.
Let $F_1,\dots,F_N, G_1,\dots,G_N\in\FC$\rom, $\phi,\psi\in
C_{\mathrm b}^\infty({\Bbb R }^N)$. Then\rom, for
$\mu$-a\rom.a\rom.\  $\gamma\in\Gamma$\rom,
\begin{align*}&S^\Gamma
(\phi(F_1,\dots,F_N),\psi(G_1,\dots,G_N))(\gamma)\\&\qquad =
\sum_{i,j=1}^N \partial _i \phi(F_1(\gamma),\dots,F_N(\gamma))\,
\partial_j \psi(G_1(\gamma),\dots,G_N(\gamma))
S^\Gamma(F_i,G_j)(\gamma).\end{align*}\end{lem}

\noindent {\it Proof}. Immediate by \eqref{iouvdc}.\quad
$\blacksquare$

\begin{lem} \label{weaesiuzghg}Let the conditions of Lemma~\rom{\ref{sdawawa}} be satisfied\rom. Then\rom,
$S^\Gamma (F,G)=0$ $\mu$-a\rom.e\rom.\ for all $F,G \in\FC$ such
that $F=0$ $\mu$-a\rom.e\rom.\end{lem}

\noindent {\it Proof}. Let $F\in\FC$, $F=0$ $\mu$-a.e. Then, for
any $r>0$, \begin{align*} 0&=\int_\Gamma\left( \int_{B(r)}
A(\gamma,x)\,\gamma(dx)\right) |F(\gamma)|\,\mu(d\gamma)\\ &=
\int_\Gamma \mu(d\gamma) \int_{B(r)}zm(dx) |F(\gamma\cup
\{x\})|.\end{align*} Hence, $F(\gamma\cup\{x\})=0$ for $\mu\otimes
m$-a.e.\ $(\gamma,x)\in\Gamma\times\R$. For any fixed
$\gamma\in\Gamma$, $\R\ni x\mapsto F(\gamma\cup\{x\})\in{\Bbb R}$
is a smooth function. Hence, for $\mu$-a.e.\ $\gamma\in\Gamma$,
$F(\gamma\cup\{x\})=0$ for all $x\in\R$, and so $S^\Gamma
(F){:=}S^\Gamma(F,F)=0$ $\mu$-a.e.\ on $\Gamma$. Using the
Cauchy--Schwarz inequality, we obtain the assertion.\quad
$\blacksquare$

%\begin{rem}\rom{In fact, one can show that, under
%the conditions of Lemma~\rom{\ref{sdawawa}}, $\mu (U)>0$ for every
%$\varnothing\ne U\subset \dd\Gamma$\rom, $U$ open. Since
%$\FC\subset C(\Gamma)$, one therefore concludes that, if two
%functions from $\FC$ coincide $\mu$-a.e., then they coincide on
%the whole $\Gamma$.}\end{rem}

\begin{prop}  \label{hgddkf} Let the conditions of
Proposition~\rom{\ref{waessedf}} be fulfilled and  let  $\phi$
also satisfy \rom{(A1).}  Then\rom, $({\cal E}_\mu^\Gamma,\FC)$ is
a pre-Dirichlet form on $L^2(\Gamma;\mu)$ \rom(i\rom.e\rom{.,} if
$({\cal E}_\mu^\Gamma,\FC)$ is closable\rom, then its closure
$({\cal E}_\mu^\Gamma,D({\cal E}_\mu^\Gamma))$ is a Dirichlet
form\rom{).}\end{prop}

\noindent{\it Proof}.  Since $\FC$ is dense in $L^2(\Gamma;\mu)$,
the assertion follows by Lemmas~\ref{kjhdrd}, \ref{weaesiuzghg}
directly from \cite[Chap.~I, Proposition~4.10]{MR92} (see also
\cite[Chap.~II, Exercise~2.7]{MR92}).\quad $\blacksquare$

\section{Closability of the (pre-)Dirichlet form ${\cal E}_\mu^\Gamma$}
 \label{ersews6543}

In this section, we shall prove that, under some condition on the
potential $\phi$, the bilinear form  $({\cal E}_\mu^\Gamma,\FC)$
is closable on $L^2(\Gamma;\mu)$. So, in what follows, we suppose
the following:

\begin{description}

\item[(A2)] Let
$$ \Phi(x){:=}\phi(x)\vee 0,\qquad x\in\R,$$ and for $\mu\otimes
m$-a.e.\ $(\gamma,x)\in\Gamma\times \R$ we set
$$\rho(\gamma,x){:=}\exp\left[-\sum_{y\in\gamma}\Phi(x-y)\right].$$
Then,  for $\mu$-a.e.\ $\gamma\in\Gamma$, $ \rho(\gamma,\cdot)=0$
$m$-a.e.\ on $$\R\setminus\big\{\,x\in\R\mid \int_{\Lambda_x}
\rho(\gamma,\cdot)^{-1}\, dm<\infty\text{ for some open
neighborhood $\Lambda_x$ of $x$}\,\big\}.$$

\end{description}

\begin{rem}\rom{Let us suppose that  $\phi\in C(\R\setminus\{0\})$. Then, condition (A2)
is  satisfied if there exists  $R>0$ such that $\phi(x)\le0$ for
$|x|\ge R$. Indeed, in this case, for each $\gamma\in\Gamma$,
$\rho(\gamma,\cdot)$ is a positive continuous function on
$\R\setminus\gamma$, which evidently yields (A2). Alternatively,
if $\mu$ is a Ruelle measure, for (A2) to hold it suffices that,
for each $\gamma\in S_\infty$, the series
$\sum_{y\in\gamma}\Phi(\cdot-y)$ converges locally uniformly on
$\R\setminus\gamma$. For Gibbs measures in low activity-high
temperature regime, in the latter condition the set $S_\infty$ can
be replaced by the set of all configurations $\gamma\in\Gamma$
satisfying, for all $N\in\N$, $$ |\gamma_{\Lambda_N}|\le C(\gamma)
 m(\Lambda_N),\qquad  C(\gamma)>0,$$ for a fixed sequence
 $\{\Lambda_N\}_{N=1}^\infty\subset {\cal O}_c(\R)$ such that $\Lambda_N\subset
 \Lambda_{N+1}$, $m(\Lambda_{N+1}\setminus\Lambda_N)\ge N+1$, $N\in\N$, and
 $\bigcup_{N=1}^\infty\Lambda_N=\R$,
which is also a set of full $\mu$ measure (cf.\
\cite[Theorem~5.2.4]{Kuna}, see also \cite[Proposition~1]{BS}).
}\end{rem}

\begin{th}\label{2323342340997656} Let the conditions of
Proposition~\rom{\ref{waessedf}} be fulfilled and  let  $\phi$
also satisfy \rom{(A1)} and \rom{(A2).} Then\rom, the bilinear
form $({\cal E}_\mu^\Gamma,\FC)$ is closable on
$L^2(\Gamma;\mu)$\rom.
\end{th}

\noindent{\it Proof}. Let $$ \Psi(x){:=}\phi(x)\wedge 0,\qquad
x\in\R .$$  We now define an auxiliary bilinear form
\begin{equation}\label{awapo}{\cal
E}^\Gamma_{\mu,\Psi}(F,G){:=}\int_\Gamma \mu(d\gamma)\int_{\R}
\gamma(dx)\,
\exp\left[\sum_{y\in\gamma\setminus\{x\}}\Psi(x-y)\right] \la
\nabla_x F(\gamma),\nabla_x G(\gamma)\ra,\end{equation} where
$F,G\in D({\cal E}^\Gamma_{\mu,\Psi}){:=}\FC$. By \eqref{waeaw},
\eqref{iouvdc}, \eqref{awapo}, and Lemma~\ref{sdawawa}
\begin{equation}\label{qwewpo}{\cal E}^\Gamma_\mu(F)\ge {\cal
E}^\Gamma_{\mu,\Psi}(F),\qquad F\in\FC,\end{equation} and
particularly the bilinear form ${\cal E}^\Gamma_{\mu,\Psi}$ is
well defined.  Using \eqref{fdrtsdrt}, \eqref{awapo}, and
Lemma~\ref{sdawawa}, we get
\begin{align}{\cal E}^\Gamma_{\mu,\Psi}(F)&=\int_\Gamma
\mu(d\gamma)\int_{\R}zm(dx)\,
\exp\left[-\sum_{y\in\gamma}\phi(x-y)+\sum_{y\in\gamma}\Psi(x-y)\right]
|\nabla_x F(\gamma+\eps_x)|^2\notag\\ &=\int_\Gamma
\mu(d\gamma)\int_{\R}z
m(dx)\,\exp\left[-\sum_{y\in\gamma}\Phi(x-y)\right] |\nabla_x
F(\gamma+\eps_x)|^2.\label{32653}\end{align}

{\it Claim}. The bilinear form $({\cal E}^\Gamma_{\mu,\Psi},\FC)$
is closable on $L^2(\Gamma;\mu)$.

For $\mu$-a.e.\ $\gamma\in\Gamma$, we define a measure
$\sigma_\gamma(dx){:=}\rho(\gamma,x)m(dx)$ on $\R$, and we
introduce the following biliear form on the space
$L^2(\R;\sigma_\gamma )$: $$ {\cal
E}_{\sigma_\gamma}(f,g){:=}\int_\R\langle\nabla f,\nabla
g\rangle\, d\sigma_\gamma,\qquad  f,g\in{\cal
D}^{\sigma_\gamma},$$ where ${\cal D}^{\sigma_\gamma}$ denote  the
$\sigma_\gamma$-classes determined by $\cal D$. Then, by
\cite[Theorem~5.3]{AR90} or \cite[Theorem~6.2]{Luis}, it follows
from (A2) that the form $({\cal E}_{\sigma_\gamma},{\cal
D}^{\sigma_\gamma})$ is closable for $\mu$-a.e.\
$\gamma\in\Gamma$. Notice also that
$$\exp\left[\sum_{y\in\gamma\setminus\{x\}}\Psi(x-y)\right]\le1
\qquad \text{for $\mu$-a.e.\ $\gamma\in\Gamma$}. $$ Now, the proof
of the claim is completely analogous to the proof of
\cite[Theorem~6.3]{Luis} (see also the survey paper
\cite{Roeckner}).

 Let $(F_n)_{n=1}^\infty$ be a sequence in
$\FC$ such that
\begin{equation}\label{qjiuz}\|F_n\|_{L^2(\mu)}\to0\quad\text{as
}n\to\infty\end{equation} and \begin{equation}\label{7654646}
{\cal E}^\Gamma_\mu(F_{n}-F_{k})\to0\quad\text{as
}n,k\to\infty.\end{equation} To prove the closability of ${\cal
E}^\Gamma_{\mu}$, it suffices to show that there exits a
subsequence $(F_{n_k})_{k=1}^\infty$ such that $$ {\cal
E}_\mu^\Gamma(F_{n_k})\to 0\quad\text{as }k\to\infty.$$ By
\eqref{qwewpo} and \eqref{7654646},
\begin{equation}\label{23856}{\cal E}^\Gamma_{\mu,\Psi}(F_n-F_k)\to
0\quad\text{as }n,k\to\infty.\end{equation} By the claim, the form
${\cal E}^\Gamma_{\mu,\Psi}$ is closable on $L^2(\Gamma;\mu)$, and
therefore \eqref{qjiuz} and \eqref{23856} imply that $${\cal
E}_\Psi(F_n)\to 0\quad\text{as }n\to\infty .$$ From here and
\eqref{32653} $$ \int_\Gamma \mu(d\gamma)\int_{\R}z
m(dx)\,\exp\left[-\sum_{y\in\gamma}\Phi(x-y)\right] |\nabla_x
F_n(\gamma+\eps_x)|^2\to0\quad \text{as } n\to\infty.$$ Therefore,
there exists a subsequence $(F_{n_k})_{k=1}^\infty$ such that
\begin{equation}\label{awawlkl} \exp\left[-\sum_{y\in\gamma}\Phi(x-y)\right] |\nabla_x
F_{n_k}(\gamma+\eps_x)|^2\to0\qquad
 \text{as $k\to\infty$ for
$\mu\otimes m$-a.e.\ } (\gamma,x)\in\Gamma\times\R.\end{equation}
By Lemma~\ref{sdawawa}, \begin{equation}
\label{awqppiiou}\sum_{y\in\gamma}|\Phi(x-y)|\le
\sum_{y\in\gamma}|\phi(x-y)|<\infty\qquad\text{for $\mu\otimes
m$-a.e.\ }(\gamma,x)\in\Gamma\times\R.\end{equation} Thus, by
\eqref{awawlkl} and \eqref{awqppiiou}, $$|\nabla_x
F_{n_k}(\gamma+\eps_x)|\to0\quad \text{as $k\to\infty$ for
$\mu\otimes m$-a.e.\ } (\gamma,x)\in\Gamma\times\R.$$ Therefore,
by Fatou's lemma, \begin{align*}{\cal
E}^\Gamma_\mu(F_{n_k})&=\int_\Gamma\mu(d\gamma)\int_{\R}zm(dx)|\nabla_x
F_{n_k}(\gamma+\eps_x)|^2\\
&=\int_\Gamma\mu(d\gamma)\int_{\R}zm(dx)|\nabla_x
F_{n_k}(\gamma+\eps_x)-\lim_{l\to\infty}\nabla_xF_{n_l}(\gamma+\eps_x)|^2\\
&\le
\liminf_{l\to\infty}\int_\Gamma\mu(d\gamma)\int_{\R}zm(dx)|\nabla_x
F_{n_k}(\gamma+\eps_x)-\nabla_xF_{n_l}(\gamma+\eps_x)|^2\\&=\liminf_{l\to\infty}
{\cal E}^\Gamma_\mu(F_{n_k}-F_{n_l}),\end{align*} which by
\eqref{7654646} can be made arbitrarily small for $k$ large
enough.\quad $\blacksquare$\vspace{2mm}

In what follows, we shall denote by $({\cal E}_\mu^\Gamma,D({\cal
E}_\mu^\Gamma))$ the closure of $({\cal E}_\mu^\Gamma,\FC)$.

\section{Another condition of closability of the form ${\cal E}_\mu^\Gamma$,
the generator of ${\cal E}_\mu^\Gamma$}\label{Section5}

Though we have given in Section~\ref{ersews6543} a condition on
$\phi$ ensuring the closability of\linebreak $({\cal
E}_\mu^\Gamma,\FC)$, we have no information about the generator of
$({\cal E}_\mu^\Gamma,D({\cal E}_\mu^\Gamma))$, except for the
fact that it exists.
 In this section, under a stronger restriction  on the
growth of the potential $\phi$ at  zero, we shall show that the
domain of the generator  contains $\FC$, and we shall give an
explicit formula for the action of the generator on this set.

Let us introduce the following condition on the potential $\phi$:

\begin{description}\item[(A3)] Let $r>0$ be as in (A1).
We have $$\int_{B(r)}e^{\phi(x)}\, m(dx)<\infty.$$
\end{description}

Notice that condition (A3) still admits that $\phi(x)\to+\infty$
as $x\to0$.

\begin{th} \label{rerere} Let the conditions of
Proposition~\rom{\ref{waessedf}} be fulfilled and  let  $\phi$
also satisfy \rom{(A1)} and \rom{(A3).} Then\rom, for any
$F,G\in\FC$\rom, $${\cal
E}^\Gamma_\mu(F,G)=\int_{\Gamma}({H}^\Gamma_\mu
F)(\gamma)G(\gamma)\,\mu(d\gamma),$$ where
\begin{gather*}({H}^\Gamma_\mu
F)(\gamma){:=}-\sum_{x\in\gamma}\exp\left[\sum_{y\in\gamma\setminus\{x\}}\phi(x-y)\right]
\Delta_x F(\gamma)\qquad \text{\rom{for $\mu$-a.e.\
$\gamma\in\Gamma$}},\\ \Delta_x F(\gamma){:=}\Delta_y
F(\gamma-\varepsilon_x+\varepsilon_y)\big|_{y=x},\end{gather*}
$\Delta$ denoting the Laplacian on $\R$\rom, and ${H}^\Gamma_\mu$
is an operator in $L^2(\Gamma;\mu)$ with domain\linebreak
$D({H}^\Gamma_\mu){:=}\FC$\rom.
\end{th}

\begin{cor}\label{deseses} Under the
conditions of Theorem~\rom{\ref{rerere}}\rom, the bilinear form
$({\cal E}_\mu^\Gamma,\FC)$ is closable on $L^2(\Gamma;\mu)$ and
the operator $({H}_\mu^\Gamma,\FC)$ has a Friedrichs
extension\rom,
 which we denote by $({H}_\mu^\Gamma,D({H}_\mu^\Gamma))$\rom.\end{cor}

\begin{rem}\rom{Notice that, in the above theorem, we do not demand condition (A2)
to hold.  }\end{rem}

\noindent {\it Proof of Theorem\/}~\ref{rerere}. First, we note
that, for any $F\in\FC$ and $\gamma\in\Gamma$, the function
$f(x){:=}F(\gamma+\eps_x)-F(\gamma)$ belongs to $\cal D$ and
$\nabla f(x)=\nabla_x F(\gamma+\eps_x)$. Therefore,  we have, for
any $F,G\in\FC$,
\begin{align*}{\cal E}^\Gamma_\mu(F,G)&=\int_\Gamma
\mu(d\gamma)\int_{\R}zm(dx)\,\la \nabla_x
F(\gamma+\eps_x),\nabla_x G(\gamma+\eps_x)\ra\\ &=-\int_\Gamma
\mu(d\gamma)\int_{\R}zm(dx)\, \Delta_x
F(\gamma+\eps_x)G(\gamma+\eps_x)\\ &=\int_\Gamma ({ H}_\mu^\Gamma
F)(\gamma)G(\gamma)\,\mu(d\gamma),\end{align*} and we have to show
that ${H}_\mu^\Gamma F\in L^2(\Gamma;\mu)$.

As easily seen, it suffices to prove that, for any
$\Lambda\in{\cal O}_c(\R)$, \begin{equation}\label{a1}
\int_\Gamma\left( \sum_{x\in\gamma_\Lambda}\exp\left[
\sum_{y\in\gamma\setminus\{x\}}\phi(x-y)\right]\right)^2
\,\mu(d\gamma)<\infty. \end{equation} By (A1), (A3) and
\eqref{gfuzag} \begin{equation}\label{ztdf56}
\int_{\R}|1-e^{\varphi(x)}|\, m(dx)<\infty.\end{equation} Hence,
using \cite[Lemma~5.2]{KK} and  Proposition~\ref{waessedf}, we get
\begin{gather} \int_\Gamma \sum_{x\in\gamma_\Lambda}\exp\left[
2\sum_{y\in\gamma\setminus\{x\}}\phi(x-y)\right]
\,\mu(d\gamma)\notag\\ =\int_\Lambda zm(dx)\int_\Gamma
\mu(d\gamma)\exp\left[\sum_{y\in\gamma}\phi(x-y)\right]\notag\\=
\int_{\Lambda}zm(dx)\left( 1+\sum_{n=1}^\infty\frac1{n!}
\int_{(\R)^n}\prod_{i=1}^n
(e^{\phi(x-y_i)}-1)k_\mu^{(n)}(y_1,\dots,y_n)\, m(dy_1)\dotsm
m(dy_n)\right)\notag\\  \le \int_\Lambda
zm(dx)\left(1+\sum_{n=1}^\infty
\frac{\xi^n}{n!}\int_{(\R)^n}\prod_{i=1}^n |1-e^{\phi(x-y_i)}|\,
m(dy_1)\dotsm m(dy_n)\right)\notag\\ =zm(\Lambda)
\sum_{n=0}^\infty
\frac{\xi^n}{n!}\left(\int_{\R}|1-e^{\phi(y)}|\,m(dy)\right)^n<\infty.\label{a2}
\end{gather}
Next, applying equality \eqref{fdrtsdrt} twise, we get from
\eqref{ztdf56}:
\begin{gather} \int_\Gamma
\sum_{x_1\in\gamma_\Lambda}\sum_{x_2\in\gamma_\Lambda\setminus\{x_1\}}\exp\left[
\sum_{y_1\in\gamma\setminus\{x_1\}}\phi(x_1-y_1)+\sum_{y_2\in\gamma\setminus
\{x_2\} }\phi(x_2-y_2)\right]\,\mu(d\gamma)\notag\\ =\int_\Gamma
\mu(d\gamma)\int_\Lambda zm(dx_1)\int_\Lambda
zm(dx_2)\exp\left[-\sum_{y_1\in\gamma}\phi(x_1-y_1)-\sum_{y_2\in\gamma\cup\{x_1\}}\phi(x_2-y_2)
\right. \notag\\ \text{}\left.\vphantom{\sum_{y_1\in\gamma}}
+\sum_{y_1\in\gamma\cup\{x_2\}}\phi(x_1-y_1)+\sum_{y_2\in\gamma\cup\{x_1\}}\phi(x_2-y_2)\right]
\notag\\ =\int_\Gamma \mu(d\gamma)\int_\Lambda
zm(dx_1)\int_\Lambda zm(dx_2)\, e^{\phi(x_1-x_2)}\notag\\ \le
z^2m(\Lambda)\int_{\R}|1-e^{\phi(x)}|\,m(dx)+z^2m(\Lambda)^2<\infty.\label{a3}\end{gather}
From \eqref{a2} and \eqref{a3} we conclude \eqref{a1}, and so the
theorem is proved.\quad $\blacksquare$

\begin{cor} \label{cfdfdp} Let the conditions of
Proposition~\rom{\ref{waessedf}} be fulfilled and  let  $\phi$
also satisfy either \rom{(A1), (A2)} or \rom{(A1), (A3).}
Then\rom, $({\cal E}_\mu^\Gamma, D({\cal E}^\Gamma_\mu))$ is a
Dirichlet form on $L^2(\Gamma;\mu)$\rom.
\end{cor}

\noindent {\it Proof}. Immediate by Proposition~\ref{hgddkf},
Theorem~\ref{2323342340997656}, and Corollary~\ref{deseses}.\quad
$\blacksquare$

\section{Quasi-regularity and diffusions}\label{Section6}
The diffusion process corresponding to the Dirichlet form $({\cal
E}_\mu^\Gamma, D({\cal E}^\Gamma_\mu))$ will, in general, live on
the bigger state space $\dd\Gamma$ consisting of all $\Z_+$-valued
Radon measures on $\R$ (which is Polish, see e.g.\ \cite{Ka75}).
Since $\Gamma\subset\dd\Gamma$ and ${\cal
B}(\dd\Gamma)\cap\Gamma={\cal B}(\Gamma)$, we can consider $\mu $
as a measure on $(\dd\Gamma,{\cal B}(\dd\Gamma))$ and
correspondingly $({\cal E}_\mu^\Gamma, D({\cal E}^\Gamma_\mu))$ as
a Dirichlet form on $L^2(\dd\Gamma;\mu)$.

The definition of quasi-regularity given in \cite[Chap.~IV,
Def.~3.1]{MR92} obviously simplifies now as follows: $({\cal
E}_\mu^\Gamma, D({\cal E}^\Gamma_\mu))$ on $L^2(\dd\Gamma;\mu)$ is
quasi-regular if and only if there exists an ${\cal
E}_\mu^\Gamma$-nest $(K_n)_{n\in\N}$ consisting of compact sets in
$\dd\Gamma$.

We recall that a sequence $(A_n)_{n\in\N}$ of closed subsets of
$\dd\Gamma$ is called an ${\cal E}_\mu^\Gamma$-nest if
$$\big\{\, F\in D({\cal E}_\mu^\Gamma)\mid F=0\text{ on
$\dd\Gamma\setminus A_n$ for some $n\in\N$}\,\big\}$$ is dense in
$D({\cal E}^\Gamma_\mu)$ with respect to the norm
$$\|\cdot\|_{{\cal E}^\Gamma_{\mu,1}}{:=}\left({\cal E
}_\mu^\Gamma(\cdot,\cdot)+(\cdot,\cdot)_{L^2(\dd\Gamma,\mu)}\right)^{1/2}.$$

\begin{prop} \label{2487743} Under the
conditions of Corollary~\rom{\ref{cfdfdp},} the Dirichlet form
$({\cal E}_\mu^\Gamma, D({\cal E}^\Gamma_\mu))$ is
quasi-regular\rom.  \end{prop}

\noindent {\it Proof}. By \cite[Proposition~4.1]{MR98}, it
suffices to show that there exists a bounded, complete
 metric $\overline\rho$ on $\dd\Gamma$ generating the vague
 topology such that, for all $\gamma\in\dd\Gamma$, $\overline\rho(\cdot,\gamma)
\in D({\cal E}_\mu^\Gamma)$ and
$S^\Gamma(\overline\rho(\cdot,\gamma))\le \eta$ $\mu$-a.e.\ for
some $\eta\in L^1(\dd\Gamma;\mu)$ (independent of $\gamma$). The
proof of this fact is quite analogous to the proof of
\cite[Proposition~4.8]{MR98}. Let us outline the main changes
needed.

We introduce the  space $\cal V$ as the completion of $\FCo$ with
respect to the norm $$ |F|_\Gamma{:=}\left(\int
S^\Gamma(F)\,d\mu\right)^{1/2}+\int |F|\,d\mu, \qquad F\in\FCo.$$
Here, $\FCo$ denotes the set of all functions on $\dd\Gamma$ of
the form \eqref3. The formulation and the proof of
\cite[Lemma~4.2]{MR98} now carries over to our case. In
particular, $\cal V$ is continuously embedded into
$L^1(\Gamma;\mu)$ and $S^\Gamma$ extends  uniquely to a bilinear
map from $({\cal V},|\cdot|_\Gamma)\times ({\cal
V},|\cdot|_\Gamma)$ into $L^1(\dd\Gamma;\mu)$.

Lemma~4.3 in \cite{MR98} now reads as follows: Let $f\in{\cal D}$.
Then, $\la f,\cdot\ra\in{\cal V}$ and $$ S^\Gamma(\la
f,\cdot\ra)(\gamma)=\int A(\gamma,x)S(f)(x)\,\gamma(dx)\qquad
\text{for $\mu$-a.e.\ $\gamma\in\dd\Gamma$}.$$ Here,
$S(f){:=}S(f,f)$ and $S(f,g)(x){:=}\la \nabla f(x),\nabla
g(x)\ra$, $f,g\in{\cal D}$, $x\in\R$. The proof is again the same.

Next, we consider the following norm on $\cal D$:
$$|f|_E{:=}\left(\int S(f)\, dm\right)^{1/2}+\int|f|k_\mu^{(1)}\,
dm,$$ where $k_\mu^{(1)}$ is the first correlation function of
$\mu$. Recall that, by Proposition~\ref{waessedf},
$k_\mu^{(1)}\le\xi$ for some $\xi>0$. Furthermore, by
\eqref{fdrtsdrt} and Lemma~\ref{sdawawa},
$$k_\mu^{(1)}(x)=\int_\Gamma
\exp\left[-\sum_{y\in\gamma}\phi(x-y)\right]\, \mu(d\gamma),\qquad
x\in\R, $$ and consequently, applying once more
Lemma~\ref{sdawawa}, we see that $k^{(1)}_\mu(x)>0$ for $m$-a.e.\
$x\in\R$. Thus, the measures $m$ and $k_\mu^{(1)}m$ are
equivalent.

 We evidently have \begin{equation}\label{ftdtdertedrdrs}|f|_E\ge |\la
f,\cdot\ra|_\Gamma\qquad \text{for all $f\in{\cal
D}$}.\end{equation} Let $\overline {\cal D}$ denote the completion
of $\cal D$ with respect to $|\cdot|_E$. The counterpart of
\cite[Lemma~4.4]{MR98} in our case reads as follows:

\begin{lem}\label{rzrtr} The inclusion map $i:({\cal
D},|\cdot|_E)\subset(L^1(\R;
k_\mu^{(1)}m),\|\cdot\|_{L^1(\R;k_\mu^{(1)}m)})$ extends uniquely
to a continuous inclusion $\overline i:\overline {\cal
D}\hookrightarrow L^1(\R;k_\mu^{(1)}m)$\rom. Furthermore\rom, $S$
extends uniquely to a bilinear continuous map from $(\overline
{\cal D},|\cdot|_E )\times (\overline {\cal D},|\cdot|_E)$ to
$L^1(\R;m)$ satisfying \rom{(S1)--(S3)} in \rom{\cite{MR98}} with
$\cal D$ replaced with $\overline {\cal D}$\rom.\end{lem}

\noindent  {\it Proof}. Let $f_n\in{\cal D}$, $n\in\N$, be an
$|\cdot|_E$-Cauchy sequence such that $f_n\to 0$ in $L^1(\R;
k_\mu^{(1)}m)$ as $n\to\infty$. Then, by \eqref{ftdtdertedrdrs},
$(\la f_n,\cdot\ra)_{n\in\N}$ is a $|\cdot|_\Gamma$-Cauchy
sequence in $\cal V$ such that $\la f_n,\cdot\ra\to 0$ in
$L^1(\Gamma;d\mu)$ as $n\to\infty$.  Hence, by what has been
proved above, $|\la f_n,\cdot\ra|_\Gamma\to 0$ as $n\to\infty$.
Therefore, $\int S^\Gamma(\la f_n,\cdot\ra)\, d\mu\to0$ as $n\to
\infty$, and so $\int S(f_n)\, dm\to0$ as $n\to\infty$.
Consequently, $|f_n|_E\to0$ as $n\to\infty$. The remaining parts
of the assertion can then be easily shown.\quad
$\blacksquare$\vspace{2mm}

We define \begin{multline*}\FCC{:=}\big\{\, g(\la
f_1,\cdot\ra,\dots,\la f_N,\cdot\ra)\mid N\in\N,\\
\text{$f_1,\dots,f_N$ $m$-versions of elements in $\overline {\cal
D}$, $g\in C_{\mathrm b}^\infty({\Bbb
R}^N)$}\,\big\}.\end{multline*}

The proof of the following assertion is the same as that of
\cite[Proposition~4.6]{MR98} if one uses Lemma~\ref{rzrtr} instead
of \cite[Lemma~4.4]{MR98}:

\begin{lem}\label{ewew} We have $\FCC\subset D({\cal E}_\mu^\Gamma)$ and for
$F=g_F(\la f_1,\cdot\ra,\dots,\la f_N,\cdot\ra)$, $G=g_G(\la
g_1,\cdot\ra,\dots,\la g_M,\cdot\ra)\in\FCC$ \begin{align*}
S^\Gamma (F,G)(\gamma)&=\sum_{i=1}^N \sum_{j=1}^M \partial_i
g_F(\la f_1,\cdot\ra,\dots,\la f_N,\cdot\ra)\\ &\quad \times
\partial_j g_G(\la g_1,\cdot\ra,\dots,\la g_M,\cdot\ra) \int
A(\gamma,x)S(f_i,g_j)(x)\,\gamma(dx)\end{align*} for
$\mu$-a\rom.e\rom.\ $\gamma\in\dd\Gamma$\rom. Furthermore\rom, for
all $f\in\overline {\cal  D} $\rom, $\la f,\cdot\ra\in{\cal V }$
and $$S^\Gamma(\la f,\cdot\ra)(\gamma)=\int
A(\gamma,x)S(f)(x)\,\gamma(dx)\qquad \text{\rom{for $\mu$-a.e.\
$\gamma\in\dd\Gamma$}}. $$\end{lem}

Next, a counterpart of \cite[Lemma~4.7]{MR98} can be easily
formulated and proved. In particular, as  norm on $\overline{\cal
D}\cap L^2(\R;k_\mu^{(1)}\,dm)$ we take $$|f|_{E,2}{:=}\left( \int
S(f)\, dm+\int f^2k_\mu^{(1)}\,dm\right)^{1/2}.$$

The formulation and the proof of \cite[Lemma~4.10]{MR98} now
remain  without essential changes, and therefore, using
Lemma~\ref{ewew}, we get analogously to the proof of
\cite[Lemma~4.11]{MR98} that $F_k\in D({\cal E}_\mu^\Gamma)$ and
\begin{equation}\label{seepwg} S^\Gamma(F_k)(\gamma)\le\int
A(\gamma,x)\tilde \chi^2_{j_k}(x)\,\gamma(dx)\quad \text{for
$\mu$-a.e.\ $\gamma\in\dd\Gamma$},\end{equation} where $F_k$ and
$\tilde\chi_{j_k}$ are as in \cite{MR98} (the function
$F_k(\cdot)=F_k(\cdot,\gamma_0)$ is defined for any fixed
$\gamma_0\in\dd\Gamma$, while the function $\tilde\chi_{j_k}$ is
independent of $\gamma_0$).

Finally, we set
$$c_k{:=}\left(1+\int\tilde\chi_{j_k}^2\,dm\right)^{-1/2}2^{-k/2},\qquad
k\in\N,$$ (since each $\tilde \chi_{j_k}$ is bounded and has
compact support, $\int \tilde\chi_{j_k}^2\,dm<\infty$). Evidently,
$c_k\to0$ as $k\to\infty$. We define $$
\overline\rho(\gamma_1,\gamma_2){:=}\sup_{k\in\N}c_k
\big(F_k(\gamma_1,\gamma_2)\big),\qquad
\gamma_1,\gamma_2\in\dd\Gamma.$$ By \cite[Theorem~3.6]{MR98},
$\overline\rho$ is a bounded, complete metric on $\dd\Gamma$
generating the vague topology.

Analogously to the proof of \cite[Proposition~4.8]{MR98}, we
conclude by \eqref{seepwg} that, for any fixed
$\gamma_0\in\dd\Gamma$, $\overline\rho(\cdot,\gamma_0)\in D({\cal
E}_\mu^\Gamma)$ and $$
S^\Gamma(\overline\rho(\cdot,\gamma_0))\le\eta\qquad
\text{$\mu$-a.e.},$$ where
$$\eta(\gamma){:=}\sup_{k\in\N}\left(2^{-k}\left(1+\int\tilde\chi_{j_k}^2\,dm\right)^{-1}
\int A(\gamma,x)\tilde\chi_{j_k}^2(x)\,\gamma(dx)\right).$$
Evidently, $$\int \eta\,d\mu\le \sum_{k=1}^\infty
2^{-k}\left(1+\int\tilde\chi_{j_k}^2\, dm\right)^{-1}\int \tilde
\chi_{j_k}^2\,dm<\infty,$$ which concludes the proof of the
proposition.\quad $\blacksquare$

\begin{prop}\label{44363} Under the
conditions of Corollary~\rom{\ref{cfdfdp},} $({\cal
E}_\mu^\Gamma,D({\cal E}_\mu^\Gamma))$ has the local property
\rom(i\rom.e\rom{.,} ${\cal E}_\mu^\Gamma(F,G)=0$ provided $F,G\in
D({\cal E}_\mu^\Gamma)$ with $\operatorname{supp}(|F|\mu)\cap
\operatorname{supp}(|G|\mu)=\varnothing$\rom{).}\end{prop}

\noindent{\it Proof}. Identical to the proof of
\cite[Proposition~4.12]{MR98}.\quad $\blacksquare$\vspace{2mm}

As a consequence of Propositions \ref{2487743}, \ref{44363} and
\cite[Chap.~IV, Theorem~3.5, and Chap.~V, Theorem~1.11]{MR92}, we
obtain the main result of this section.

 \begin{th}\label{8435476} Let the conditions of
Proposition~\rom{\ref{waessedf}} be fulfilled and  let $\phi$\rom,
in addition\rom, satisfy either \rom{(A1), (A2)} or \rom{(A1),
(A3).} Then\rom, there exists a conservative diffusion process
\rom(i\rom.e\rom{.,} a conservative strong Markov process with
continuous sample paths\rom) $${\bf M}=({\pmb{ \Omega}},{\bf
F},({\bf F}_t)_{t\ge0},({\pmb \Theta}_t)_{t\ge0}, ({\bf
X}(t))_{t\ge 0},({\bf P }_\gamma)_{\gamma\in\dd\Gamma})$$ on
$\dd\Gamma$ \rom(cf\rom.\ \rom{\cite{Dy65})} which is properly
associated with $({\cal E}_\mu^\Gamma,D({\cal
E}_\mu^\Gamma))$\rom, i\rom.e\rom{.,} for all \rom($\mu$-versions
of\/\rom) $F\in L^2(\dd\Gamma;\mu)$ and all $t>0$ the function
\begin{equation}\label{zrd9665} \dd\Gamma\ni\gamma\mapsto
p_tF(\gamma){:=}\int_{\pmb\Omega} F({\bf X}(t))\, d{\bf
P}_\gamma\end{equation} is an ${\cal
E}_\mu^\Gamma$-quasi-continuous version of
$\exp(-t{H}_\mu^\Gamma)F$\rom, where $H_\mu^\Gamma$ is the
generator of $({\cal E}_\mu^\Gamma,D({\cal E}_\mu^\Gamma))$
\rom(cf\rom.\ \rom{\cite[Chap.~1, Sect.~2]{MR92}).} $\bf M$ is up
to $\mu$-equivalence unique \rom(cf\rom.\ \rom{\cite[Chap.~IV,
Sect.~6]{MR92}).} In particular\rom, $\bf M$ is $\mu$-symmetric
\rom(i\rom.e\rom{.,} $\int G\, p_tF\, d\mu=\int F \, p_t G\, d\mu$
for all $F,G:\dd\Gamma\to{\Bbb R}_+$\rom, ${\cal
B}(\dd\Gamma)$-measurable\rom) and has $\mu$ as an invariant
measure\rom.
\end{th}

In the above theorem, $\bf M$ can be taken to be canonical, i.e.,
$\pmb\Omega=C([0,\infty)\to\dd\Gamma)$, ${\bf
X}(t)(\omega){:=}\omega(t)$, $t\ge 0$, $\omega\in\pmb\Omega$,
$({\bf F}_t)_{t\ge 0}$ together with $\bf F$ is the corresponding
minimum completed admissible family (cf.\
\cite[Section~4.1]{Fu80}) and ${\pmb \Theta}_t$, $t\ge0$, are the
corresponding natural time shifts.

We recall that by  $(H_\mu^\Gamma, D(H_\mu^\Gamma))$ we denote the
generator of the closed form $({\cal E}_\mu^\Gamma,D({\cal
E}_\mu^\Gamma))$.

\begin{th}\label{fdresrear} $\bf M$ from
Theorem~\rom{\ref{8435476}} is  up to $\mu$-equivalence
\rom(cf\rom.\ \rom{\cite[Definition~6.3]{MR92}}\rom) unique
between all diffusion processes ${\bf M}'=({\pmb{ \Omega}}',{\bf
F}',({\bf F}'_t)_{t\ge0},({\pmb \Theta}'_t)_{t\ge0}, ({\bf
X}'(t))_{t\ge 0},({\bf P }'_\gamma)_{\gamma\in\dd\Gamma})$ on $\dd
\Gamma$ having $\mu$ as an invariant measure and solving the
martingale problem for $(-H_\mu^\Gamma, D(H_\mu^\Gamma))$\rom,
i\rom.e\rom.\rom, there exists a set $\Gamma_0\in{\cal
B}(\dd\Gamma)$ such that $\dd\Gamma\setminus\Gamma_0$ is ${\cal
E}_\mu^\Gamma$-exceptional \rom(so\rom, in particular\rom,
$\mu(\Gamma_0)=1$\rom) and such that for all $G\in
D(H_\mu^\Gamma)$ $$\widetilde G({\bf X}'(t))-\widetilde G({\bf
X}'(0))+\int_0^t (H_\mu^\Gamma G)({\bf X}'(s))\,ds,\qquad t\ge0,$$
is an $({\bf F}_t')$-martingale under ${\bf P}_\gamma'$ for all
$\gamma\in\Gamma_0$\rom. \rom(Here\rom, $\widetilde G$ denotes a
quasi-continuous version of $G$\rom, cf\rom. \rom{\cite[Ch.~IV,
Proposition~3.3]{MR92}.)}
\end{th}

\noindent {\it Proof}. The statement of the theorem follows
directly from (the proof of) \cite[Theorem~3.5]{AR}.\quad
$\blacksquare$\vspace{2mm}

Our next aim is to show that the diffusion  process $\bf M$
properly associated with $({\cal E}_\mu^\Gamma,D({\cal
E}_\mu^\Gamma))$ lives, in fact, on the  space
$\Gamma=\Gamma_{{\Bbb R}^d}$ provided $d\ge 2$.

\begin{th}\label{12345} Let the conditions of Theorem~\rom{\ref{8435476}}
be satisfied and let $d\ge 2$\rom. Then the set
$\dd\Gamma\setminus \Gamma$ is ${\cal
E}_\mu^\Gamma$-exceptional\rom. \end{th}

\noindent {\it Proof}. We modify the proof of \cite[Proposition~1 and
Corollary~1]{RS98} according to our situation. For the convenience of
the reader we shall present the proof completely.

It suffices to prove the result locally, that is to show that, for every
fixed $a\in\N$,  the set
$$N{:=}\big\{\,\gamma\in\dd\Gamma: \sup(\gamma(\{x\}): x\in
[-a,a]^d)\ge 2\,\big\}$$ is ${\cal E}_\mu^\Gamma$-exceptional. By
\cite[Lemma~1]{RS98}, we need to prove that there exists a
sequence $(u_n)_{n=1}^\infty\subset D({\cal E}_\mu^\Gamma)$ such
that each $u_n$, $n\in\N$, is a continuous function on
$\dd\Gamma$, $u_n\to {\bf 1} _{N}$ pointwise as $n\to\infty$, and
$\sup_{n\in\N}{\cal E}_\mu^\Gamma (u_n)<\infty$.

Let $f\in C_0^\infty ({\Bbb R})$ be such that ${\bf 1}_{[0,1]}\le
f\le {\bf 1}_{[-1/2,3/2)}$ and $|f'|\le 3 \times{\bf 1}_{[-1/2,3/2)}$. For
any $n\in\N$ and $i=(i_1,\dots,i_d)\in\Z^d$, define a function $f_i^{(n)}
\in{\cal D}$ by $$ f_i^{(n)}(x){:=}\prod_{k=1}^d f(n x_k- i_k),\qquad
  x\in\R.$$ Let also $I_i^{(n)}(x){:=}\prod_{k=1}^d {\bf 1}
_{[-1/2,3/2)}(nx_k-i_k)$, $x\in\R$, and note that $f_i^{(n)}\le
I_i^{(n)}$.  Since $$\partial _j f_i^{(n)}(x)=n f' (nx_j-i_j)
\prod_{k\ne j} f(nx_k-i_k),$$ we get \begin{equation}\label{one}
|\nabla f_i^{(n)}(x)|^2\le 9n^2 d I_i^{(n)}(x).\end{equation}

Let $\psi\in C_{\mathrm b}^\infty ({\Bbb R})$ be such that ${\bf
1} _{[2,\infty)}\le\psi\le{\bf 1}_{[1,\infty)}$ and $|\psi'|\le
2\times {\bf 1}_{ (1,\infty)}$. Set ${\cal
A}{:=}\Z^d\cap[-na,na]^d$ and define continuous functions
$$\dd\Gamma\ni\gamma\mapsto u_n(\gamma){:=}\psi \left(
\sup_{i\in{\cal A}}\la f_i^{(n)},\gamma\ra\right),\qquad n\in\N.$$
Evidently, $u_n\to {\bf 1}_N$ pointwise as $n\to\infty$.
Furthermore, by an appropriate approximation of the function
$${\Bbb R}^{|{\cal A}|} \ni (y_1,\dots,y_{|{\cal A}|})\mapsto
\sup_{i\in{\cal A}}y_i$$ by $C_{\mathrm b}^\infty ({\Bbb
R}^{|{\cal A}|})$ functions  (compare with \cite[Lemma~3.2]{RS95}
and \cite[Lemma~4.7]{MR98}), we conclude that, for each $n\in\N$,
$u_n\in{D}({\cal E}_\mu^\Gamma)$ and
\begin{equation}\label{two}S^\Gamma (u_n)(\gamma)\le \left(\psi'\left(
\sup_{i\in{\cal A}}\la
f_i^{(n)},\gamma\ra\right)\right)^2\sum_{x\in\gamma}
A(\gamma,x)\sup_{i\in{\cal A}}|\nabla
f_i^{(n)}(x)|^2\qquad\text{for $\mu$-a.e.\ $\gamma\in\dd\Gamma$}.
\end{equation} Next, we have for each $\gamma\in\dd\Gamma$
\begin{equation} \left(\psi'\left(
\sup_{i\in{\cal A}}\la f_i^{(n)},\gamma\ra\right)\right)^2
                    \le 4 \times{\bf 1}_{\{\sup_{i\in{\cal A}}
                    \la f_i^{(n)},
\gamma\ra>1\}}\le 4\times{\bf 1}_{\{\sup_{i\in{\cal A}}\la
I_i^{(n)},\gamma\ra\ge2\}},\label{three}\end{equation} where we
used the fact that $\la I_i^{(n)},\gamma\ra$ is an integer. Thus,
by \eqref{one}--\eqref{three} \begin{align} S^\Gamma (u_n)(\gamma)
&\le 4 \times {\bf 1}_{\{\sup_{i\in{\cal A}}\la
I_i^{(n)},\gamma\ra\ge2\}}
\sum_{x\in\gamma}A(\gamma,x)\sup_{i\in{\cal A}} \big(9n^2d
I_i^{(n)}(x)\big)\notag\\ &\le 36 n^2 d\sum_{i\in{\cal A}} {\bf
1}_{\{\la I_i^{(n)},\gamma\ra\ge2\}}  \sum_{x\in\gamma}
A(\gamma,x){\bf 1}_{[-a-1,a+1] ^d}(x)\qquad\text{for $\mu$-a.e.\
$\gamma\in\dd\Gamma$ }.\notag\end{align} Consequently,
\begin{align}\int S^\Gamma(u_n)\,d\mu&\le 36 n^2 d\sum_{i\in{\cal A}}\int
_{\dd\Gamma}\mu(d\gamma)\int_{[-a-1,a+1]^d}zm(dx)\,{\bf 1}_{\{ \la
I_i^{(n)},\gamma+\eps_x\ra\ge 2\}}(\gamma,x)\notag\\ &= 36 n^2
d\sum_{i\in{\cal A}}\int_{\dd\Gamma}\mu(d\gamma)\bigg(
\int_{\{I_i^{(n)}=1\}}zm(dx)\, {\bf 1}_{\{\la
I_i^{(n)},\gamma\ra\ge 1\}}(\gamma)
\notag\\&\quad\qquad+\int_{[-a-1,a+1]^d\setminus\{I_i^{(n)}=1\}}zm(dx)\,
{\bf 1}_{\{\la I_i^{(n)},\gamma\ra\ge2\}}(\gamma)\bigg)\notag\\
&\le 36n^2 dz \sum_{i\in{\cal A}}\big(m(\{I_i^{(n)}=1\})\mu (\{\la
I_i^{(n)},\cdot\ra\ge1\})\notag\\&\quad\qquad +(2a+2)^d \mu(\{\la
I_i^{(n)},\cdot\ra\ge2\})\big).\label{four}\end{align}

By using \cite[Theorem~5.5]{Ru70}, we easily conclude that there exist
constants $c_1,c_2\in(0,\infty)$, independent of $i$ and $n$, such that
for all $i\in{\cal A}$ and $n\in\N$
\begin{align}\mu(\{\la I_i^{(n)},\cdot\ra\ge1\})&\le c_1 m(\{I_i^{(n)}=1\}),
\notag\\
\mu(\{\la I_i^{(n)},\cdot\ra\ge2\})&\le c_2 m(\{I_i^{(n)}=1\})^2
.\label{five}
\end{align} Thus, by \eqref{four} and \eqref{five}, there exists $c_3\in
(0,\infty)$, independent of $n$, such that for all $n\in\N$
 \begin{equation}\label{six}
\int S^\Gamma(u_n)\,d\mu \le c_3 n^2 \sum_{i\in{\cal
A}}m(\{I_i^{(n)}=1\})^2.
\end{equation}
Since $|{\cal A}|=(2na+1)^d$ and $m(\{I_i^{(n)}=1\})=(2/n)^d$, we finally
get from \eqref{six}: \begin{align*} {\cal E}_\mu^\Gamma(u_n)&\le c_3 n^2
(2na+1)^d \bigg(\frac 2n\bigg)^{2d}\\&\le \operatorname{const}_d\qquad
\text{for all }n\in\N\end{align*} for some $\operatorname{const}_d\in(0,
\infty)$, provided $d\ge2$. \quad $\blacksquare$\vspace{2mm}

As a direct consequence of Theorem~\ref{12345}, we get

\begin{cor}\label{rews}  Let the conditions of Theorem~\rom{\ref{8435476}}
be satisfied and let $d\ge 2$\rom. Then\rom, the assertions of
Theorems~\rom{\ref{8435476}} and \rom{\ref{fdresrear}} hold with
$\dd \Gamma$ replaced by $\Gamma$\rom. \end{cor}

\section{Scaling limit of the stochastic dynamics}\label{Section7}

Throughout this section, we shall suppose that  $\phi$ satisfies
(S) and  (UI) with $z=1$ and $\mu\in{\cal G}(1,\phi)$ is the
measure corresponding to the construction with  empty boundary
conditions. We shall now discuss a scaling limit of the diffusion
process constructed in Theorem~\ref{8435476}, the scaling being
absolutely analogous to the one considered in
\cite{Brox,Rost,Spohn,GP,GKLR}.

\subsection{Scaling of the process}

First, let us briefly recall a result of Brox \cite{Brox} on a
scaling limit of Gibbs measures.

{\it First scaling}. We scale the position of the particles inside
the configuration space $\Gamma$ as follows:
$$\Gamma\ni\gamma\mapsto S_{{\mathrm
in},\,\epsilon}(\gamma){:=}\{\epsilon x\mid
x\in\gamma\}\in\Gamma,\qquad \epsilon>0.$$ Let us define the image
measure
\begin{equation}\label{cdsres}\tilde\mu_\epsilon{:=}S_{{\mathrm
in},\,\epsilon}^*\mu.\end{equation} As easily seen,
$\tilde\mu_\epsilon$ is an element of ${\cal
G}(\epsilon^{-d},\phi_\epsilon)$ with
$\phi_\epsilon{:=}\phi(\epsilon^{-1}\cdot)$. Furthermore,
$\tilde\mu_\epsilon$ satisfies (UI) and corresponds to the
construction with empty boundary conditions.

{\it Second scaling}. This scaling leads us out of  the
configuration space and is given by $$\Gamma\ni \gamma\mapsto\So
(\gamma){:=}\epsilon^{d/2}(\gamma-k^{(1)}_{\tilde\mu_\epsilon
}m)\in\Gamma_\epsilon,\qquad \epsilon>0, $$ where
$\Gamma_\epsilon{:=} \So(\Gamma)\subset{\cal D}'$, and as before
${\cal D}'$ is the topological dual of $\cal D$ (where both $\cal
D$ and ${\cal D}'$ are equipped with their respective usual
locally convex topologies). We consider $\Gamma_\epsilon$ as a
topological subspace of ${\cal D}'$, thus $\Gamma_\epsilon$ is
equipped with the corresponding Borel $\sigma$-algebra. Obviously,
$\So:\Gamma\to\Gamma_\epsilon$ is continuous, hence
Borel-measurable. Since it is also one-to-one and since both
$\Gamma$ and ${\cal D}'$ are standard measurable spaces, it
follows by \cite[Chap.~V, Theorem~2.4]{Par67} that
$\Gamma_\epsilon$ is a Borel subset of ${\cal D}'$ and that
$\So^{-1}:\Gamma_\epsilon\to\Gamma$ is also Borel-measurable. The
function $k^{(1)}_{\tilde\mu_\epsilon}$ is the first correlation
function of $\tilde\mu_\epsilon$.  It easily follows from
\eqref{cdsres} that
$$k^{(1)}_{\tilde\mu_\epsilon}=\epsilon^{-d}\rho,$$ where
$\rho{:=}k_\mu^{(1)}$ is the first correlation function of the
measure $\mu$ (which is a constant because of the translation
invariance of the measure $\mu$). Thus,
$$\So(\gamma)=\epsilon^{d/2}\gamma- \epsilon^{-d/2}\rho
m{=:}\gamma_\epsilon.$$ We now set
\begin{equation}\label{dfrdeess}
\mu_\epsilon{:=}\So^*\tilde\mu_\epsilon=\So^*S_{{\mathrm
in},\,\epsilon}^*\mu.\end{equation}

Let $$u_\mu^{(2)}(x_1,x_2){:=}k_\mu^{(2)}(x_1,x_2)-
k_\mu^{(1)}(x_1)k_\mu^{(2)}(x_2)=k_\mu^{(2)}(x_1,x_2)-\rho^2 $$
denote the second Ursell function of the measure $\mu$. By
\cite[Theorem~4.5]{Brox} or \cite[Chapter~4]{Ru69}, we have
$$\int_\R |u_\mu^{(2)}(x,0)|\,m(dx)<\infty,$$ and let
$$c{:=}\rho+\int_\R u_\mu^{(2)}(x,0)\,m(dx)$$ (which is the
compressibility of the Gibbs state $\mu$). We define a Gaussian
measure $\nu_c$ on $({\cal D }',{\cal B}({\cal D}'))$ by its
Fourier transform $$\int_{{\cal
D}'}\exp\big(i\la\varphi,\omega\ra\big)\,\nu_c(d\omega)=\left(-\frac
c2\,\int_\R\varphi(x)^2\,m(dx)\right),\qquad \varphi\in{\cal D}.$$
We have the following result (cf.\ \cite[Theorem~6.5]{Brox}):

\begin{prop} \label{wejktfjhg} Let us assume that the potential $\phi$
satisfies \rom{(S)} and \rom{(UI)} with $z=1$\rom, and let
$\mu\in{\cal G }(1,\phi)$  be the Gibbs measure corresponding to
the construction with empty boundary conditions\rom. For each
$\epsilon>0$\rom, consider $\mu_\epsilon$\rom, defined by
\eqref{dfrdeess}\rom,  as a probability measure on $({\cal
D}',{\cal B}({\cal D}'))$\rom. Then\rom, the family of measures
$(\mu_\epsilon)_{\epsilon>0}$ converges weakly on ${\cal D}'$ to
the Gaussian measure $\nu_c$\rom. \end{prop}

For simplicity of notation, in what follows  we shall exclude the
case $d=1$. However, all our further considerations do also work
in that case.

The scaled process of our interest is defined by $$ {\bf
X}_\epsilon(t){:=}S_{{\mathrm out},\,\epsilon}(S_{{\mathrm
in},\,\epsilon}{\bf X}(\epsilon^{-2}t)),\qquad t\ge0,\
\epsilon>0,$$ where $({\bf X}(t))_{t\ge0}$ is the process
constructed in Corollary~\ref{rews}. Next, for each $\epsilon>0$,
we construct a Dirichlet form ${\cal E}_\epsilon$ such that $({\bf
X}_\epsilon(t))_{t\ge0}$ is the unique process which is properly
associated to ${\cal E}_{\epsilon}$.

Since the transformation $S_{{\mathrm out},\,\epsilon}$ is
invertible,  we can define a unitary operator ${\cal S}_{{\mathrm
out },\,\epsilon}: L^2(\Gamma_\epsilon;\mu_\epsilon)\to L^2
(\Gamma;\tilde\mu_\epsilon)$ by setting ${\cal S}_{{\mathrm out
},\,\epsilon} F$ to be the $\tilde\mu_\epsilon$-class represented
by $\widetilde F\circ \So$ for any $\mu_\epsilon$-version
$\widetilde F$ of $F\in L^2(\Gamma_\epsilon;\mu_\epsilon)$. Using
this operator, we define a bilinear form $({\cal
E}_\epsilon,D({\cal E }_\epsilon))$ on
$L^2(\Gamma_\epsilon,\mu_\epsilon)$ as the image of the bilinear
form $({\cal E}^\Gamma_{\tilde\mu_\epsilon},D({\cal
E}^\Gamma_{\tilde\mu_\epsilon}))$ under ${\cal S}_{{\mathrm out
},\,\epsilon}^{-1}$: \begin{equation}\label{wawawawa}{\cal
E}_\epsilon (F,G){:=}{\cal E}^\Gamma_{\tilde\mu_\epsilon} ({\cal
S}_{{\mathrm out },\,\epsilon} F,{\cal S}_{{\mathrm out
},\,\epsilon} G),\qquad F,G\in D({\cal E}_\epsilon),\end{equation}
where $ D({\cal E}_\epsilon){:=} {\cal S}_{{\mathrm out
},\,\epsilon}^{-1}\big(D({\cal
E}^\Gamma_{\tilde\mu_\epsilon})\big) $. It follows from
\cite[Chapter~VI, Exercise~1.1]{MR92} that $({\cal
E}_\epsilon,D({\cal E}_\epsilon))$ is a Dirichlet form.
 Let ${\cal H}_\epsilon$ (respectively $\widetilde{\cal H}_\epsilon$)
 denote the generator of
the form $({\cal E}_\epsilon,D({\cal E}_\epsilon))$ (respectively
$({\cal E}^\Gamma_{\tilde\mu_\epsilon},D({\cal
E}^\Gamma_{\tilde\mu_\epsilon}))$) on $L^2(\Gamma;\mu_\epsilon)$
(respectively $L^2(\Gamma;\tilde\mu_\epsilon)$).  Then, it follows
from the definition of $({\cal E}_\epsilon,D({\cal E}_\epsilon))$
that \begin{equation}\label{sewaswawap}{\cal H}_\epsilon={\cal
S}_{{\mathrm out },\,\epsilon}^{-1}\widetilde{\cal H}_\epsilon
{\cal S}_{{\mathrm out },\,\epsilon}.\end{equation}

We have the following proposition (compare with
\cite[Theorem~4.1]{GKLR}).

\begin{prop}\label{poseski} Let the potential $\phi$
fulfill  conditions \rom{(S), (UI)} with $z=1$, \rom{(A1)} and
either \rom{(A2)} or \rom{(A3)} and let $\mu\in{\cal G}(1,\phi)$
be the Gibbs measure constructed as a limit of finite volume Gibbs
measures  with empty boundary conditions\rom. For
$\omega\in\Gamma_\epsilon$\rom, let ${\bf
Q}^\epsilon_\omega{:=}{\bf P}_{S_{{\mathrm
in},\,\epsilon}^{-1}\So^{-1}\omega}$\rom. Then\rom, for all
\rom($\mu_\epsilon$-versions\rom) of $F\in
L^2(\Gamma_\epsilon;\mu_\epsilon)$ and all $t>0$\rom, the function
$$\Gamma_\epsilon\ni \omega\mapsto p_\epsilon(t,F)(\omega){:=}\int
_{\pmb{\Omega}}F({\bf X}_\epsilon(t))\, d{\bf
Q}^\epsilon_{\omega}$$ is a $\mu_\epsilon$-version of
$\exp(-t{\cal H}_\epsilon)F$\rom. The process \begin{equation}
\label{process} {\bf M}_\epsilon{:=}(\pmb\Omega,{\bf F},({\bf
F}_{t/\epsilon^2})_{t\ge0}, (\pmb\Theta_{t/\epsilon^2})_{t\ge0},
({\bf X}_\epsilon)(t))_{t\ge0},({\bf
Q}^\epsilon_\omega)_{\omega\in\Gamma_\epsilon})\end{equation} is a
diffusion process and thus up to $\mu_\epsilon$-equivalence the
unique process in this class which is properly associated with
$({\cal E}_\epsilon,D({\cal E}_\epsilon))$\rom.  It has $\mu_\eps$
as an invariant measure. \end{prop}

\noindent {\it Proof}. By \eqref{sewaswawap}, to prove the first
statement of the theorem it suffices to show that, for all
($\tilde\mu_\epsilon$-versions of) $F\in
L^2(\Gamma;\tilde\mu_\epsilon)$ and all $t>0$\rom, the function
$$\Gamma\ni \gamma\mapsto \int _{\pmb{\Omega}}F(S_{{\mathrm in
},\,\epsilon}({\bf X}(\epsilon^{-2}t)))\, d{\bf P}_{S_{{\mathrm
in},\,\epsilon}^{-1}\gamma},\qquad $$ is a
$\tilde\mu_\epsilon$-version of $\exp(-t\widetilde{\cal
H}_\epsilon)F$\rom. Analogously to ${\cal S}_{{\mathrm out},\,
\epsilon}$, we define a unitary operator ${\cal S}_{{\mathrm
in},\, \epsilon}: L^2(\Gamma;\tilde\mu_\epsilon)\to
L^2(\Gamma;\mu)$ through the transformation $S_{{\mathrm in},\,
\epsilon}$ of $\Gamma$. We  note that $\FC$ remains invariant
under ${\cal S}_{{\mathrm in},\, \epsilon}$. A direct calculation
shows that, for any $F,G\in\FC$,
\begin{equation}\label{fdes}{\cal E}_\mu^\Gamma ({\cal S}_{{\mathrm in},\, \epsilon}F,
{\cal S}_{{\mathrm in},\, \epsilon}G)=\epsilon^2 {\cal
E}_{\tilde\mu_\epsilon}^\Gamma(F,G).\end{equation} Since ${\cal
S}_{{\mathrm in},\, \epsilon}^{-1}H_\mu^\Gamma {\cal S}_{{\mathrm
in},\, \epsilon}$ is the generator of the closure of the bilinear
form \linebreak $({\cal E}_\mu^\Gamma({\cal S}_{{\mathrm in},\,
\epsilon}\cdot,{\cal S}_{{\mathrm in},\, \epsilon}\cdot),\FC)$ on
$L^2(\gamma,\tilde\mu_\epsilon)$, \eqref{fdes} implies that $$
{\cal S}_{{\mathrm in},\, \epsilon}^{-1}H_\mu^\Gamma{\cal
S}_{{\mathrm in},\ \epsilon}=\epsilon^2
H^\Gamma_{\tilde\mu_\epsilon}.$$ By using Theorem~\ref{8435476},
we now easily obtain the first assertion.

The fact that ${\bf M}_\epsilon$ is a diffusion is straightforward
to check.  In particular, it then follows from \cite[Chap.~IV,
Theorem~3.5]{MR92} that ${\bf M}_\epsilon$ is properly associated
with $({\cal E}_\epsilon,D({\cal E}_\epsilon))$. \quad
$\blacksquare$\vspace{2mm}

\subsection{Scaling limit of the Dirichlet form ${\cal E}_\mu^\Gamma$}

We shall now show  the convergence of the processes ${\bf
M}_\epsilon$ to a generalized Ornstein--Uhlenbeck process in the
sense of convergence of the corresponding Dirichlet forms ${\cal E
}_\epsilon$. The limiting Dirichlet form will coincide, up to a
constant factor, with the limiting Dirichlet form of \cite{GKLR}.

We introduce the set $\FCD$ of all   functions on ${\cal D }'$ of
the form \eqref3 where $\Gamma$ is replaced by ${\cal D}'$. Thus,
any function $F\in\FC$ is a restriction of some $\widetilde
F\in\FCD$ to $\Gamma$. Notice that any function from  ${\cal
S}_{{\mathrm out },\,\epsilon}^{-1}\big(\FC\big)$, defined on
$\Gamma_\epsilon$, may be extended to a function from $\FCD$, and
so the set (of $\mu_\epsilon$-classes of) $\FCD$ is dense in
$D({\cal E}_\epsilon)$ with respect to the norm $\|\cdot\|_{{\cal
E }_\epsilon}{:=}\big(\|\cdot\|_{L^2(\mu_\epsilon)}^2+{\cal E
}_\epsilon(\cdot)\big)^{1/2}$.

We next introduce a bilinear form ${\cal E}_{\nu_c}$ on $L^2({\cal
D}';\nu_c)$ as follows:\begin{align*}{\cal
E}_{\nu_c}(F,G)&=\int_{{\cal D}'}\int_{\R}\partial _x
F(\omega)(-\Delta_x)\partial_xG(\omega)\, m(dx)\nu_c(d\omega)\\
&=\int_{{\cal D}'}\int_{\R}\la\nabla_x\partial_x
F(\omega),\nabla_x\partial_x G(\omega)\,
m(dx)\,\nu_c(d\omega),\end{align*} where $F,G\in D({\cal
E}_{\nu_c}){:=}\FCD$. Here, $\partial _x$ denotes the derivative
in direction $\eps_x$, i.e., $$\partial _x F(\omega)=\frac d{dt}\,
F(\omega+t\eps_x)\big|_{t=0}, $$ and $\nabla_x$ and $\Delta_x$
denote the gradient  and the Laplacian in the $x$ variable,
respectively. One easily sees that, for $F\in\FCD$ of the form
\eqref{3},
\begin{multline}\label{55}{\cal E}_{\nu_c}(F){:=}{\cal
E}_{\nu_c}(F,F)=\sum_{i,j=1}^N\int_{{\cal D}'}  \partial_i
g_F(\la\varphi_1,\omega\ra,\dots,\la \varphi_N,\omega\ra)
\\ \times\partial_j g_F(\la\varphi_1,\omega\ra,\dots,\la
\varphi_N,\omega\ra)\,\nu_c(d\omega) \int_\R\la
\nabla\varphi_i(x),\nabla\varphi_j(x)\ra\, m(dx).\end{multline}

By using the integration by parts formula on Gaussian space (e.g.\
\cite[Ch.~6, Theorems~6.1.2 and 6.1.3]{BK}), we conclude that
$${\cal E}_{\nu_c}(F,G)=\int_{{\cal D}'}({\cal
H}_{\nu_c}F)(\omega) G(\omega)\,\nu_c(d\omega),$$ where
\begin{align*}({\cal H}_{\nu_c}F)(\omega):&= -\sum_{i,j=1}^N
\partial_i\partial_j g_F(\la\varphi_1,\omega\ra,\dots,\la \varphi_N,\omega\ra)
\int_\R \la\nabla\varphi_i(x),\nabla\varphi_j(x)\ra\, m(dx)\\
&\quad-c^{-1}\sum_{j=1}^N\partial_j
g_F(\la\varphi_1,\omega\ra,\dots,\la \varphi_N,\omega\ra)
\la\Delta\varphi_j,\omega\ra.\end{align*} Hence, the bilinear form
${\cal E}_{\nu_c}$ is closable on $L^2({\cal D}',\nu_c)$.
Moreover, it is well known (e.g.\ \cite[Ch.~6, Theorem~6.1.4]{BK})
that the operator  ${ \cal H}_{\nu_c}$ is essentially self-adjoint
on $\FCD$. We preserve the same notation  for its closure. The
operator ${ \cal H}_{\nu_c}$ generates an infinite-dimensional
Ornstein--Uhlenbeck semigroup $(\exp(-t{ \cal
H}_{\nu_c}))_{t\ge0}$ in $L^2({\cal D}',\nu_c)$. This semigroup is
associated to a generalized Ornstein--Uhlenbeck process $({\bf
N}(t))_{t\ge0}$ on ${\cal D}'$, see e.g. \cite[Chapter~6,
Section~1.5]{BK}. This process informally satisfies the stochastic
differential equation \eqref{uhisidf}.

\begin{th} \label{wawardersuih} Let the conditions of
Proposition~\rom{\ref{poseski}} be fulfilled\rom. Then\rom, the
bilinear forms ${\cal E}_\epsilon$ converge to the bilinear form
${\cal E}_{\nu_c}$ in the following sense\rom: for all
$F,G\in\FCD$\rom, $${\cal E}_\eps (F,G)\to{\cal
E}_{\nu_c}(F,G)\quad \text{\rom{as} }\epsilon\to0.$$
\end{th}

\begin{rem}\rom{Proposition~\ref{poseski} and Theorem~\ref{wawardersuih}  tell us that
the scaled process $({\bf X}_\epsilon(t))_{t\ge0}$ converges to
the Ornstein--Uhlenbeck process $({\bf N}(t))_{t\ge0}$ in the
sense of the convergence of their respective Dirichlet forms on
$\FCD$. In particular, equation~\eqref{uhisidf} allows us  to
identify $c^{-1}$ as the bulk diffusion coefficient corresponding
to the initial process $({\bf X}(t))_{t\ge0}$.}
\end{rem}

\noindent{\it Proof of Theorem\/} \ref{wawardersuih}. Due to the
polarization identity, it suffices to show that, for each
$F\in\FCD$, $${\cal E}_\epsilon(F)\to{\cal
E}_{\nu_c}(F)\quad\text{as }\epsilon\to0.$$

Let $F\in\FCD$ be of the form \eqref{3}. Then, by
\eqref{wawawawa},
\begin{gather}{\cal E}_\epsilon (F)=
\int_\Gamma \tilde\mu_\epsilon (d\gamma)\int_\R \epsilon ^{-d}
m(dx)\, |\nabla_x ({\cal S}_{{\mathrm
out},\epsilon}F)(\gamma+\eps_x)|^2\notag\\ =\int_\Gamma
\tilde\mu_\epsilon (d\gamma)\int_\R \epsilon ^{-d} m(dx)\,\big|
\nabla_x \big(g_F(\la
\varphi_1,\epsilon^{d/2}(\gamma+\eps_x)-\rho\epsilon^{-d/2}\ra,\notag\\
\dots,\la \varphi_N
,\epsilon^{d/2}(\gamma+\eps_x)-\rho\epsilon^{-d/2}\ra)\big)\big|^2\notag\\
=\int_\Gamma \tilde\mu_\epsilon (d\gamma)\int_\R \epsilon ^{-d}
m(dx)\,  \sum_{i,j=1}^N \partial_i g_F(\la
\varphi_1,\epsilon^{d/2}\gamma-\rho\epsilon^{-d/2}\ra+\epsilon^{d/2}\varphi_1(x),\notag\\
\dots,\la
\varphi_N,\epsilon^{d/2}\gamma-\rho\epsilon^{-d/2}\ra+\epsilon^{d/2}\varphi_N(x))
\,\partial_j g_F(\la
\varphi_1,\epsilon^{d/2}\gamma-\rho\epsilon^{-d/2}\ra+\epsilon^{d/2}\varphi_1(x),\notag\\
\dots,\la
\varphi_N,\epsilon^{d/2}\gamma-\rho\epsilon^{-d/2}\ra+\epsilon^{d/2}\varphi_N(x))
\la
\epsilon^{d/2}\nabla\varphi_i(x),\epsilon^{d/2}\nabla\varphi_j(x)\ra
\notag\\ =\sum_{i,j=1}^N\int_{{\cal D}'}\mu_\epsilon
(d\omega)\int_\R m(dx)\,
\partial_i g_F
(\la\varphi_1,\omega\ra+\epsilon^{d/2}\varphi_1(x),
\dots,\la\varphi_N,\omega\ra+\epsilon^{d/2} \varphi_N(x))\notag \\
\times \partial_j g_F
(\la\varphi_1,\omega\ra+\epsilon^{d/2}\varphi_1(x),
\dots,\la\varphi_N,\omega\ra+\epsilon^{d/2} \varphi_N(x))\la
\nabla\varphi_i(x),\nabla\varphi_j(x)\ra.
\label{534493}\end{gather}

Let $\hat\mu_\epsilon$, resp.\ $\hat\nu_c$ denote the measure on
${\Bbb R}^N$ obtained as the image of $\mu_\epsilon$, resp.\
$\nu_c$ under the mapping $${\cal D}'\ni\omega\mapsto (\la
\varphi_1,\omega\ra,\dots,\la \varphi_N,\omega\ra)\in{\Bbb R}^N.$$
Then, it follows  from Proposition~\ref{wejktfjhg} that $\hat
\mu_\epsilon$ converges weakly on ${\Bbb R}^N$ to $\hat\nu_c$.
Since the functions $\partial_i g_F$, $i=1,\dots,N$, are
continuous and bounded on ${\Bbb R}^N$, we therefore get from
\eqref{55}:
\begin{multline}\label{77}\sum_{i,j=1}^N \int_{{\Bbb R}^N}
\partial_i g_F(x_1,\dots,x_N)\,\partial_j g_F(x_1,\dots,x_N)\,
\hat\mu_\epsilon(dx_1,\dots,dx_N)\\ \times \int_\R \la
\nabla\varphi_i(x),\nabla\varphi_j(x)\ra\, m(dx)\to {\cal
E}_{\nu_c}(F)\quad\text{as }\epsilon\to0.\end{multline} Choose
$\Lambda\in{\cal O}_c(\R)$  such that
$\operatorname{supp}\varphi_i\subset\Lambda$ for all
$i=1,\dots,N$. Then,  by \eqref{534493} and \eqref{77}, it
suffices to show that, for any $i,j\in\{1,\dots,N\}$,
\begin{multline}\label{71} \int_{{\Bbb R}^N}\hat\mu_\epsilon
(dx_1,\dots,dx_N)\int_\Lambda m(dx)\, \big|
\partial_i
g_F(x_1+\epsilon^{d/2}\varphi_1(x),\dots,x_N+\epsilon^{d/2}\varphi_N(x))\\
\times
\partial_j
g_F(x_1+\epsilon^{d/2}\varphi_1(x),\dots,x_N+\epsilon^{d/2}\varphi_N(x))\\
\text{}
-
\partial_i g_F(x_1,\dots,x_N)\,\partial_j g_F(x_1,\dots,x_N)
  \big|\to0\quad
\text{as }\epsilon\to0.\end{multline}

Set $\Lambda_n{:=}(-n,n)^N$, $n\in\N$. Since $\hat\nu_c({\Bbb R}
^N)=1$, we get $\hat\nu_c(\Lambda_n)\to1$ as $n\to\infty$. Hence,
for any fixed $\delta>0$, there exists $n_0\in\N$ such that
$\hat\nu_c(\Lambda_{n_0})\ge 1-\delta$. Since $\Lambda_{n_0}$ is
open and $\hat\mu_\epsilon$ converges weakly to $\hat\nu_c$, we
conclude from the Portemanteau theorem that
$$\liminf_{\epsilon\to0}\hat\mu_\epsilon (\Lambda_{n_0})\ge
\hat\nu_c(\Lambda_{n_0})\ge1-\delta.$$ Hence, there exists
$\epsilon_0>0$ such that, for each $\epsilon<\epsilon_0$,
$\hat\mu_\epsilon(\Lambda_{n_0})\ge 1-2\delta$, so that $\hat
\mu_\epsilon (\Lambda_{n_0}^c)\le 2\delta$. From here
\begin{multline}\label{88} \int_{\Lambda_{n_0}^c}\hat\mu_\epsilon
(dx_1,\dots,dx_N)\int_\Lambda m(dx)\, \big|
\partial_i
g_F(x_1+\epsilon^{d/2}\varphi_1(x),\dots,x_N+\epsilon^{d/2}\varphi_N(x))\\
\times
\partial_j
g_F(x_1+\epsilon^{d/2}\varphi_1(x),\dots,x_N+\epsilon^{d/2}\varphi_N(x))
- \partial_i g_F(x_1,\dots,x_N)\,\partial_j g_F(x_1,\dots,x_N)
  \big|\\ \le  4\delta m(\Lambda)\max_{i=1,\dots,N}\sup_{(x_1,\dots,x_N)\in{\Bbb
  R}^N}|\partial_i g_F(x_1,\dots,x_N)|^2\end{multline} for
  $\epsilon<\epsilon_0$.

Let $\alpha{:=}\max_{i=1,\dots,N}\sup_{x\in\R}|\varphi_i(x)|$.
Since the function $\partial_i g_F\,\partial_j g_F$ is uniformly
continuous on the compact set $[-n_0-\alpha,n_0+\alpha]^N$, we
conclude that there exists $\epsilon_1>0$, $\epsilon_1\le
\epsilon_0$, such that \begin{multline*} \big|
\partial_i
g_F(x_1+\epsilon^{d/2}\varphi_1(x),\dots,x_N+\epsilon^{d/2}\varphi_N(x))
\,\partial_j
g_F(x_1+\epsilon^{d/2}\varphi_1(x),\dots,x_N+\epsilon^{d/2}\varphi_N(x))
\\\text{}- \partial_i g_F(x_1,\dots,x_N)\,\partial_j g_F(x_1,\dots,x_N)
  \big|<\delta \end{multline*} for all $(x_1,\dots,x_N)\in\Lambda_{n_0}$
and $\epsilon<\epsilon_1$. Therefore,
\begin{multline}\label{99}
\int_{\Lambda_{n_0}}\hat\mu_\epsilon(dx_1,\dots,dx_N)\int_\Lambda
m(dx) \, \big|
\partial_i
g_F(x_1+\epsilon^{d/2}\varphi_1(x),\dots,x_N+\epsilon^{d/2}\varphi_N(x))\\
\times
\partial_j
g_F(x_1+\epsilon^{d/2}\varphi_1(x),\dots,x_N+\epsilon^{d/2}\varphi_N(x))\\
\text{}- \partial_i g_F(x_1,\dots,x_N)\,\partial_j
g_F(x_1,\dots,x_N)
  \big|\le m(\Lambda)\delta\end{multline} for each
  $\epsilon<\epsilon_1$. Finally, \eqref{88} and \eqref{99} imply
  \eqref{71}.\quad $\blacksquare$

\subsection{Tightness} In this subsection, we shall discuss  the
problem of convergence in law of the processes ${\bf M}_\epsilon$
as $\epsilon\to0$.

For $\epsilon>0$ the law of the scaled equilibrium process is the
probability measure on $C([0,\infty),\Gamma_\epsilon)$ given by $$
{\bf P}_\epsilon {:=} {\bf Q}_{\mu_\epsilon}\circ {\bf
X}_\epsilon^{-1}, $$ where $${\bf
Q}_{\mu_\epsilon}{:=}\int_{\Gamma_\epsilon} {\bf
Q}^\epsilon_\omega \, \mu_\epsilon(d\omega),$$ (cf.\
Proposition~\ref{poseski}).  Since $C([0,\infty),\Gamma_\epsilon)$
is a Borel subset of $C([0,\infty),{\cal D}')$ (under the natural
embedding) with compatible measurable structure, we can consider
${\bf P}_\epsilon$ as a measure on the (common for all
$\epsilon>0$) space $C([0,\infty),{\cal D}')$.

For $n\in\Z$, we define a weighted Sobolev space ${\cal H}_n$ as
the closure of ${\cal D}$ with respect to the Hilbert norm
$$\|f\|_n^2=\langle f,f\rangle_n{:=} \int_{\R}A^n f(x)f(x)\,
m(dx),\qquad f\in{\cal D},$$ where $$ Af(x){:=}-\Delta f(x)+|x|^2
f(x),\qquad x\in\R.$$ We identify ${\cal H}_0= L^2(\R;m)$ with its
dual and obtain $$ {\cal D}\subset S(\R)\subset{\cal H}_n\subset
L^2(\R;m) \subset {\cal H}_n\subset S'(\R)\subset{\cal D}',\qquad
n\in\N.$$ Here, as usual $S'(\R)$ denotes the space of tempered
distributions which is the the topological dual of $S(\R)$, the
Schwartz space of smooth functions on $\R$ decaying faster than
any polynomial. Of course, ${\cal H}_{-n}$ is the topological dual
of ${\cal H}_n$ with respect to ${\cal H}_0$. For each $n\in\Z$,
the embedding ${\cal H}_n\hookrightarrow {\cal H}_{n-d}$ is of
Hilbert--Schmidt type.

\begin{th}\label{7523} Let the conditions of Proposition~\rom{\ref{poseski}}
be satisfied\rom. Then\rom, there exists $k\in\N$, $k\ge d+1$\rom,
such that the family of probability measures $({\bf
P}_\epsilon)_{\epsilon>0}$ can be restricted to the space
$C([0,\infty),{\cal H}_{-k})$. Furthermore\rom, $({\bf
P}_\epsilon)_{\epsilon>0}$ is tight on $C([0,\infty),{\cal
H}_{-k})$\rom.\end{th}

\noindent{\it Proof}. The proof of this theorem is analogous  to
the proof of \cite[Theorem~6.1]{GKLR}.

Consider the diffusion process ${\bf M}_\epsilon$, $\epsilon>0$,
on the state space ${\cal D}'$. Considering its distribution on
$C([0,\infty),{\cal D}')$, we may regard its canonical realization
\eqref{process}. So, in particular $\pmb\Omega=C([0,\infty),{\cal
D}')$, ${\bf X}(t)(\omega)=\omega(t)$, $t\ge0$, $\pmb\Theta
_t(\omega)=\omega(t+\cdot)$, and ${\bf
P}_\epsilon=\int_{\Gamma_\epsilon}{\bf Q}^\epsilon_\omega\,
\mu_\epsilon(d\omega)$. Fix $T>0$. Below, we canonically project
the process onto $\pmb\Omega_T{:=}C([0,T],{\cal D}')$ without
expressing this explicitly. We define the time reversal
$r_T(\omega){:=}\omega(T-\cdot)$.

Let $f\in{\cal D}$. It is easy to show that $\langle
f,\cdot\rangle\in D({\cal E}_\epsilon)$. By the Lyons--Zheng
decomposition, cf.\ \cite{LZ88,FOT94, LZ94}, we have, for all
$0\le t\le T$: $$\langle f, {\bf X}(t)\rangle-\langle f,{\bf
X}(0)\rangle =\frac12\, {\bf M}_t(\epsilon,f)+\frac12\, \big({\bf
M}_{T-t}(\epsilon,f)(r_T)-{\bf M}_T(\epsilon,f)(r_T)\big),$$ ${\bf
P}_\epsilon$-a.e., where $({\bf M}_t(\epsilon,f))_{0\le t\le T}$
 is a continuous $({\bf P}_\epsilon,({\bf F}_{t/\epsilon^2})_{0\le t\le
 T})$-martingale and \linebreak $({\bf M}_t(\epsilon,f)(r_T))_{0\le t\le T}$
 is a continuous $({\bf P}_\epsilon,(r_T^{-1}({\bf F}_{t/\epsilon^2}))_{0\le t\le
 T})$-martingale.
Moreover, by \eqref{534493}, $$\langle {\bf
M}(\epsilon,f)\rangle_t= 2t \int_{\R}|\nabla f(x)|^2\, m(dx), $$
as e.g.\ follows from     \cite[Theorem~5.2.3 and
Theorem~5.1.3(i)]{FOT94} (see also a remark in the proof of
\cite[Theorem~6.1]{GKLR}). Hence, by the Burkholder--Davies--Gundy
inequality and since ${\bf P}_\epsilon\circ r_T={\bf P}
_{\epsilon}$, we can find $C>0$ such that, for all $f\in{\cal D}$,
$0\le s\le t\le T$,
\begin{align}& {\Bbb E}_{{\bf P}_\epsilon}\big[|\langle f,{\bf X}(t)\rangle-\langle f,{\bf X}(s)
\rangle|^4\big]\notag\\ &\qquad \le{\Bbb E}_{{\bf P}_\epsilon}
\big[ | {\bf M }_t(\epsilon,f)-{\bf M}_s(\epsilon,f)|^4\big]
+{\Bbb E}_{{\bf P}_\epsilon} \big[ | {\bf M
}_{T-t}(\epsilon,f)(r_T)-{\bf
M}_{T-s}(\epsilon,f)(r_T)|^4\big]\notag\\ &\qquad \le C(t-s)^2
\|\,|\nabla f |\,\|_0^4.\label{te4ews342}
\end{align}

Now, we can use \eqref{te4ews342} to define $\langle f,{\bf
X}(t)\rangle-\langle f,{\bf X}(s) \rangle$ for $f\in S(\R)$ via an
approximation  as an element of $L^4(\pmb\Omega,{\bf
P}_\epsilon)$. Then, the estimate \eqref{te4ews342} holds true for
$f\in S(\R)$.

We can choose $\alpha>0$ and $k\in\N$ large enough, so that
\begin{equation}\label{76r7r} \forall f\in S(\R):\qquad \|\,|\nabla f|\,\|_0^4\le
\alpha\|f\|_{k-2d}^4.\end{equation}
 Let  $(e_i)_{i=0}^\infty$ be the sequence of Hermite
functions forming an orthonormal basis of  ${\cal H}_{k-2d}$. For
$i\in\N$, let $a_i$ denote the eigenvalue of the operator $A$
belonging to the eigenvector $e_i$. Then,
$(a_i^{k-d}e_i)_{i\in\N}$ forms an orthonormal basis in  ${\cal
H}_{-k}.$ Hence, by \eqref{te4ews342} and \eqref{76r7r},
\begin{align} &\big( {\Bbb E}_{{\bf P}_\epsilon}\big[\| {\bf
X}(t)-{\bf X}(s)\|_{-k}^4\big]\big)^{1/2}\notag\\&\qquad
=\bigg({\Bbb E}_{{\bf P }_\epsilon}\bigg[\bigg( \sum_{i=0}^\infty
a_i^{2k-2d}\big( (e_i,{\bf X}(t))_{-k}-(e_i,{\bf
X}(s))_{-k}\big)^2\bigg)^2\bigg] \bigg)^{1/2}\notag\\&\qquad=
\bigg( {\Bbb E}_{{\bf P}_\epsilon}\bigg[\bigg( \sum_{i=0}^\infty
a_i^{-2d}\big( \langle e_i,{\bf X}(t) \rangle - \langle e_i,{\bf
X}(s)\rangle \big)^2\bigg)^2\bigg]\bigg)^{1/2}\notag\\&\qquad \le
\bigg(\sum_{i=0}^\infty
a_i^{-2d}\bigg)^{1/2}\bigg(\sum_{i=0}^\infty a_i^{-2d}\, {\Bbb
E}_{{\bf P}_\epsilon}\left[ (\langle e_i,{\bf X}(t)\rangle-\langle
e_i,{\bf X}(s)\rangle )^4\right]\bigg)^{1/2}\notag\\ &\qquad\le
C'(t-s),\label{tfrt}
\end{align}
where the constant $C'{:=}(\alpha C)^{1/2}\sum_{i=0}^\infty
a_i^{-2d}$ is finite, since $A^{-d}$ is a Hilbert--Schmidt
operator.

Since, by Proposition~\ref{wejktfjhg}, $\mu_\epsilon\to\nu_c$ as
$\epsilon\to0$, now the tightness of $({\bf
P}_\epsilon)_{\epsilon>0}$  on $C([0,\infty),{\cal H}_{-k})$
follows by standard  arguments.\quad$\blacksquare$\vspace{2mm}

It follows from Theorem~\ref{7523} that there exists at least one
accumulation point $\widetilde{\bf P}$ of $({\bf
P}_\epsilon)_{\epsilon>0}$ on $C([0,\infty),{\cal H}_{-k})$, i.e.,
${\bf P}_{\epsilon_n}\to\widetilde{\bf P}$ weakly for some
subsequence $\epsilon_n\to0$. However, it is still an open
question whether the measures ${\bf P}_{\epsilon}$ converge to the
law ${\bf P}$ of the Ornstein--Uhlenbeck process $({\bf
N}(t))_{t\ge0}$, i.e., whether the measure $\widetilde{\bf P}$
must always coincide with ${\bf P}$.

\begin{center}

{\bf Acknowledgments}
\end{center}

\noindent We would like to thank Martin Grothaus for useful
discussions. Financial support of the BiBoS-Research Center,  the
DFG-Forschergruppe 399 ``Spectral Analysis, Asymptotic
Distributions and Stochastic Dynamics,'' and the SFB 611
``Singular Phenomena and Scaling in Mathematical Models'' is
gratefully acknowledged.

\noindent Yu.~Kondratiev, Fakult\"at f\"ur Mathematik,
Universit\"at Bielefeld, Postfach 10 01 31, D-33501 Bielefeld,
Germany;  Institute of Mathematics, Kiev, Ukraine; BiBoS, Univ.\
Bielefeld, Germany\\
\texttt{kondrat@mathematik.uni-bielefeld.de}\\[2mm]
 E. Lytvynov, Institut f\"{u}r Angewandte Mathematik,
Universit\"{a}t Bonn, Wegelerstr.~6, D-53115 Bonn, Germany; BiBoS,
Univ.\ Bielefeld, Germany\\
  \texttt{lytvynov@wiener.iam.uni-bonn.de}\\[2mm]
 M. R\"ockner, Fakult\"at f\"ur Mathematik, Universit\"at
Bielefeld, Postfach 10 01 31, D-33501 Bielefeld, Germany; BiBoS,
Univ.\ Bielefeld, Germany\\
\texttt{roeckner@mathematik.uni-bielefeld.de}

\end{document}